%% file: KB.tex
\documentclass[a4paper,12pt]{article}
\usepackage{xspace}
\usepackage{epic}
\usepackage[standard,thmmarks]{ntheorem}
\usepackage{placeins} 
\usepackage{caption2}
\renewcommand{\captionlabeldelim}{.}
\renewcommand{\captionfont}{\sffamily}
\renewcommand{\captionlabelfont}{\bfseries}

\makeatletter
\renewcommand{\subsection}{%
\setcounter{subsection}{\value{theorem}}%
\stepcounter{theorem}%
\@startsection
{subsection}%
{1}
{0em}
{-\baselineskip}
{-4pt}
{\bfseries\normalsize}}
\renewcommand{\p@enumi}{\thetheorem.}
\renewcommand{\p@enumii}{\thetheorem.\theenumi.}
\makeatother
\newcommand{\Rev}[1]{\ensuremath{Rev(#1)}}
\newcommand{\Weld}[1]{\ensuremath{Weld(#1)}}
\newcommand{\Rules}{\ensuremath{Rules}\xspace}
\newcommand{\SLtwo}{\ensuremath{SL2}}
\newcommand{\NRed}[1]{\ensuremath{Rble_N(#1)}}
\newcommand{\Red}[1]{\ensuremath{Rble_D(#1)}}
\newcommand{\NRevLHS}[1]{\ensuremath{Rev_N({LHS}(#1))}}
\newcommand{\DRev}[1]{\ensuremath{Rev_D(#1)}}
\newcommand{\DRevLHS}[1]{\ensuremath{Rev_D({LHS}(#1))}}
\newcommand{\WDiff}{\textit{WDiff}\xspace}
\newcommand{\RR}{\ensuremath{\mathsf{Aut}}\xspace}
\newcommand{\Store}{\ensuremath{\mathsf{S}}\xspace}
\newcommand{\Storen}[1]{\ensuremath{\mathsf{S}(#1)}}
\newcommand{\Considered}{\ensuremath{\mathsf{Considered}}\xspace}
\newcommand{\This}{\ensuremath{\mathsf{Now}}\xspace}

\newcommand{\New}{\ensuremath{\mathsf{New}}\xspace}
\newcommand{\Delete}{\ensuremath{\mathsf{Delete}}\xspace}

\newcommand{\Permanent}{\ensuremath{\mathsf{P}}\xspace}
\newcommand{\SetOfRules}[1]{\ensuremath{\mathsf{Set}(#1)}}
\newcommand{\LHS}[1]{\ensuremath{\mathsf{LHS}(#1)}}
\newcommand{\needed}{\textsf{needed}\xspace}
\newcommand\inverse{\iota}%

\newcounter{cond}

\newcommand{\uast}{^{\textstyle \ast}} 
\newcommand{\Astar}{A\uast}
\newcommand{\naturals}{\mathbb N}

\newcommand{\shl}{short-lex\xspace}
\title{Knuth--Bendix for groups\\with infinitely many rules%
\thanks{Funded by EPSRC grant no. GR/K 76597}%
}
\author{David Epstein \and Paul Sanders}
\begin{document}
\maketitle


\begin{abstract}
We introduce a new class of groups with solvable word problem,
namely groups specified by a confluent set of \shl-reducing
Knuth--Bendix rules which form a regular language. 
This simultaneously generalizes \shl-automatic groups and
groups with a finite confluent set of \shl-reducing rules.
We describe a computer program which looks for such a set of rules in
an arbitrary finitely presented group. Our main theorem is that our
computer program finds the set of rules, if it exists, given enough
time and space. (This is an optimistic description of our result---for
the more pessimistic details, see the body of the paper.)

The set of rules is embodied in a finite state automaton in two
variables. A central feature of our program is an operation,
which we call \textit{welding}, used to combine existing rules with
new rules as they are found. Welding can be defined on arbitrary
finite state automata, and we investigate this operation in abstract,
proving that it can be considered as a process which takes as input
one regular language and outputs another regular language.

In our programs we need to convert several
non-deterministic finite state automata to deterministic versions
accepting the same language.
We show how to improve somewhat on the standard subset construction,
due to special features in our case. We axiomatize these special
features, in the hope that these improvements can be used in other
applications.

The Knuth--Bendix process normally spends most of its time in
reduction, so its efficiency depends on doing reduction quickly.
Standard data structures for doing this can become very large,
ultimately limiting the set of presentations of groups which can be so
analyzed.
We are able to give a method for rapid reduction using our much
smaller two variable automaton, encoding the (usually infinite)
regular language of rules found so far.
Time taken for reduction in a given group is a small constant
times the time taken for reduction in the best schemes known
(see \cite{Holt:KBMAG}), which is not too bad since we are reducing
with respect to an infinite set of rules, whereas known schemes use a
finite set of rules.

We hope that the method described here might lead to the computation
of automatic structures in groups for which this is currently infeasible.
\end{abstract}

\input{Sec1.tex}
\input{Sec2.tex}
\input{Sec3.tex}
\input{Sec4.tex}
\input{Sec5.tex}
\input{Sec6.tex}
\input{Sec7.tex}
\input{Sec8.tex}
\input{Sec9.tex}
\bibliography{all} 
\bibliographystyle{plain}

\noindent D.B.A. Epstein\\
Mathematics Institute, University of Warwick\\ Coventry CV4 7AL, UK\\
\verb+dbae@maths.warwick.ac.uk+

\vspace{3mm}

\noindent P.J. Sanders\\
Mathematics Institute, University of Warwick\\ Coventry CV4 7AL, UK\\
\verb+pjs@maths.warwick.ac.uk+

\end{document}

%% file: Sec1.tex
\section*{Contents} To help readers find their way around the
inevitably complex structure of this paper,
we start with a brief description of each section.

\noindent {\bf 1. Introduction.} This briefly sets some of
the background for the paper and describes the motivation for this work.

\noindent{\bf 2. Our class of groups in context.}
We define the class of groups to which this paper is devoted and prove
various relations with related classes of groups.
Groups in our class satisfy our main theorem (\thref{main}), which
states that if the set of minimal \shl reducing  rules is regular,
then our program succeeds in finding the finite state automaton which
accepts these rules.

\noindent{\bf 3. Welding.} Here we describe one of the main new ideas in
this paper, namely welding. This process can be applied to any finite
state automaton. In our case it is the tool which enables us perform
the apparently impossible task of generating an infinite set of
Knuth--Bendix rules from a finite set. Welding has good properties
from the abstract language point of view (see \thref{welding unique}).
Welding has some important features.
Firstly, if an automaton starts by
accepting only pairs $(u,v)$ such that $\bar u = \bar v$ in $G$, then
the same is true after welding.
Secondly, the welded automaton can encode infinitely many distinct
equalities, even if the original only encoded a finite number.
Thirdly, the welded automaton is usually much smaller than the
original automaton.
At the end of this section we show that any
group determined by a regular set of rules is finitely presented.

\noindent{\bf 4. Standard Knuth--Bendix.} In this section,
we describe the standard Knuth--Bendix process for
string rewriting, in the form in which it is normally used to analyze
finitely presented groups and monoids. We need this as a background against
which to describe our modifications.

\noindent{\bf 5. Our version of Knuth--Bendix.}
We give a description of our Knuth--Bendix procedure.
We describe critical pair analysis, minimization of a rule and give
some brief details of our method of reduction using a two-variable
automaton which encodes the rules.

\noindent{\bf 6. Correctness of our Knuth--Bendix Procedure.}
We prove that our Knuth--Bendix procedure does
what we want it to do. The proof is not at all easy. In part the
difficulty arises from the fact that we have to not only find new
rules, but also delete unwanted rules, the latter in the interests of
computational efficiency, or, indeed, computational feasibility.
Our main tool is the concept of a Thue path (see \thref{Thue path}).
Although it is hardly possible that this is a new concept, we have not
seen elsewhere its systematic use to understand the progress of
Knuth--Bendix with time.
One hazard in programming Knuth--Bendix is that some clever manoeuvre
changes the Thue equivalence relation.
The key result here is \thref{maintain congruence}, which carefully
analyzes the effect of various operations on Thue equivalence. In fact
it provides more precise control, enabling other hazards, such as
continual deletion and re-insertion of the same rule, to be avoided.
It is also the most important step in proving our main result,
\thref{main}. This says that if our program is applied to a group
defined by a regular set of minimal \shl rules, then, given sufficient time
and space, a finite state automaton accepting exactly these rules will
eventually be constructed by our program, after which it will loop
indefinitely, reproducing the same finite state automaton (but
requiring a steadily increasing amount of space for redundant information).

\noindent{\bf 7. Fast reduction.}
We describe a surprisingly pleasant aspect
of our data structures and procedures, namely that reduction with
respect to our probably infinite set of rules can be carried out very
rapidly. Given a reducible word $w$, we can find a rule $(\lambda,\rho)$, such
that $w$ contains $\lambda$ as a subword, in a time which is linear in
the length of $w$.
Fast algorithms in computer science are often achieved by using finite
state automata, and the current situation is an example.
We explain how to construct the necessary automata and why they work.

\noindent{\bf 8. A modified determinization algorithm.}
Here we describe a modification of the standard algorithm, to be
found in every book about computing algorithms, that
determinizes a non-deterministic finite state automaton. Our version
saves space as compared with the standard one. It is well suited
to our special situation. We give axioms which enable one to see when
this improved algorithm can be used.

\noindent{\bf 9. Miscellaneous details.}
A number of miscellaneous points are discussed.
In particular, we compare our approach to that taken in
\textit{kbmag} (see \cite{Holt:KBMAG}).

\section{Introduction}

We give some background to our paper, and describe the class of groups
of interest to us here.

A celebrated result of Novikov and Boone asserts that the word problem for
finitely presented groups is, in general, unsolvable. This means that
a finite presentation of a group is known and has been written down
explicitly, with the
property that there is no algorithm whose input is a word in the
generators, and whose output states whether or not the word is
trivial. Given a presentation of a group for which one is unable to
solve the word problem, can any help at all be given by a computer?

The answer is that some help \textit{can} be given with the kind of
presentation that arises naturally in the work of many mathematicians,
even though one can formally prove that there is no procedure that
will \textit{always} help.

There are two general techniques for trying to determine, with the help
of a computer, whether two words in a group are equal or not. One is
the Todd--Coxeter coset enumeration process and the other is the
Knuth--Bendix process. Todd-Coxeter is more adapted to finite groups
which are not too large. In this paper, we are motivated by groups which
arise in the study of low dimensional topology. In particular they are
usually infinite groups, and the number of words of length $n$ rises
exponentially with $n$. For this reason, Todd--Coxeter is not much use
in practice. Well before Todd--Coxeter has had time to work
out the structure of a large enough neighbourhood of the identity in the
Cayley graph to be
helpful, the computer is out of space.

On the other hand, the Knuth--Bendix process is much better adapted to this
task, and it has been used quite extensively, particularly by Sims,
for example in connection with computer investigations into problems
related to the Burnside problem.
It has also been used to good effect by Holt and Rees in their automated
searching for isomorphisms and homomorphisms between two given
finitely presented groups (see \cite{HoltRees:Software}).
In connection with searching for a \shl-automatic structure on a group,
Holt was the first person to realize that the
Knuth--Bendix process might be the right direction to choose (see
\cite {EpsteinHoltRees}).  Knuth--Bendix will run
for ever on even the most innocuous hyperbolic triangle groups, which are
perfectly easy to understand. Holt's successful plan was to use
Knuth--Bendix for a certain amount of time, decided heuristically, and
then to interrupt Knuth--Bendix and make a guess as to the automatic
structure. One then uses axiom-checking,
a part of automatic group theory (see
\cite[Chapter 6]{WPiG}), to see whether the guess is correct.
If it isn't correct, the checking process will produce suggestions as to how
to improve the guess.
Thus, using the concept of an automatic group as a mechanism for
bringing Knuth--Bendix to a halt has been one of the philosophical bases
for the work done at
Warwick in this field almost from the beginning.
In addition to the works already cited in this paragraph, the reader
may wish to look at \cite{HoltRees:Software} and \cite{Holt:WarwickSoftware}.

For a \shl-automatic group, a minimal set of
Knuth--Bendix rules may be infinite, but it is always a regular
language (see \thref{automatic implies regular}),
and therefore can be encoded by a finite state machine.
In this paper, we carry this philosophical approach further,
attempting to compute this finite state machine directly, and to carry
out as much of the Knuth--Bendix process as possible using only
approximations to this machine.

Thus, we describe a setup that can handle an infinite regular set of
Knuth--Bendix rewrite rules. For our setup to be effective, we need to
make several assumptions. Most important is the assumption that we
are dealing with a group, rather than with a monoid. Secondly, our
procedures are perhaps unlikely to be of much help unless the group actually
is \shl-automatic.
Our main theorem---see \thref{main}---is that our
Knuth--Bendix procedure succeeds in constructing the finite state
machine which accepts the (unique) confluent set of \shl minimal rules
describing a group, if and only if this set of rules is a regular
language.

Previous computer implementations of the semi-decision procedure to
find the \shl-automatic structure on a group are
essentially specializations of the Knuth--Bendix procedure
\cite{KnuthBendix} to a
string rewriting context together with fast, but space-consuming,
automaton-based methods of performing word reduction relative to a finite set
of \shl-reducing rewrite rules. Since \shl-automaticity
of a given finite presentation is, in general, undecidable, space-efficient
approaches to the Knuth--Bendix procedure are desirable.
Our new algorithm performs a Knuth--Bendix type procedure relative
to a possibly infinite regular set of \shl-reducing rewrite rules,
together with a companion word reduction algorithm which has been designed with
space considerations in mind.

In standard Knuth--Bendix, there is a tension between time and space
when reducing words. Looking for a left-hand side in a word can take a
long time, unless the left-hand sides are carefully arranged in a data
structure that traditionally
takes a lot of space. Our technique can do very rapid
reduction without using an inordinate amount of space (although, for
other reasons, we have not been able to save as much space as we
originally hoped).
This is explained in \thref{A modified determinization algorithm}.

We would like to thank Derek Holt for many conversations about this
project, both in general and in detail. His help has, as always, been generous
and useful.

%% file: Sec2.tex
\section{Our class of groups in context}
In this paper we study groups, together with a finite ordered set of
monoid generators,  with the property that their set of
universally minimal \shl rules is a regular language.
In this section, we explain what this rather daunting sentence means,
and we set this class of groups in the context of various other related
classes, investigating which of these classes is included in which.
In the next section, we will prove that groups in this class are finitely
presented.

Throughout we will work with a group $G$ generated by a fixed finite set $A$,
and a fixed finite set of defining relations.
Formally, we are given a map $A\to G$, but our language will
sometimes (falsely) pretend that $A$ is a subset of $G$. The reader is
urged to remain aware of the distinction, remembering that, as a result of
the insolubility of the word problem, it is not in general possible to tell
whether the given map $A\to G$ is injective. We assume we are
given an involution
$\inverse:A\to A$ such that, for each $x\in A$, $\inverse(x)$ represents
$x^{-1}\in G$. By $\Astar$ we mean the set of words (strings) over
$A$. (Formally a word is a function $\{1,\ldots,n\}\to A$, where $n\ge 0$.)
We also write
$\inverse:\Astar\to\Astar$ for the formal inverse map defined by
$\inverse(x_1\ldots x_p) =\inverse(x_p)\ldots\inverse(x_1)$.

We assume we are given a fixed total order on $A$.
This allows us to define the \textit{\shl} order on $\Astar$ as
follows.
We denote by $|u|$ the length of $u \in \Astar$.
If $u,v\in \Astar$, we
say that $u<v$ if either
$|u| < |v|$ or $u$ and $v$ have the same
length and $u$ comes before $v$ in lexicographical order.
The \textit{\shl representative}
of $g\in G$ is the smallest $u\in\Astar$ such
that $u$ represents $g$. This is also called the \textit{\shl normal
form} of $g$. If $u\in\Astar$, we write $\overline u \in G$ for the element
of $G$ which it represents. If $u$ is the \shl representative of
$\overline u$, we say that $u$ is \textit{in \shl normal form}.

Suppose we have $(G,A)$ as above. Then there may or may not be an algorithm
that has a word $u\in\Astar$ as input and as output the \shl
representative of $\overline u\in G$. The existence of such an algorithm is
equivalent to the solubility of the word problem for $G$,
since there are only a finite number of words $v$ such that $v<u$.

A natural attempt to construct such an algorithm is to find a set
$R$ of \textit{replacement rules}, also known as \textit{Knuth--Bendix rules}.
In this paper, a replacement rule
will be called simply a \textit{rule}, and we will restrict our
attention to rules of a rather special kind.
A \textit{rule} is a pair $(u,v)$ with $u>v$
Given a rule $(u,v)$, $u$ is called the \textit{left-hand
side} and $v$ the \textit{right-hand side}.
The idea of the
algorithm is to start with an arbitrary word $w$ over $A$ and to \textit{reduce}
it as follows: we change it
to a smaller word by looking in $w$ for some left-hand side $u$ of some rule
$(u,v)$ in $R$. We then replace $u$ by $v$ in $w$ (this is called an
\textit{elementary reduction}) and repeat the operation until no
further elementary reductions are possible (the repeated process is
called a \textit{reduction}).
Eventually the process must stop with an
\textit{$R$-irreducible} word, that is a word which contains no subword
which is a left-hand side of $R$.

\subsection{Thue equivalence.}\label{Thue}
Given a set of rules $R$, we write $u\rightarrow_R v$ if there is an elementary
reduction from $u$ to $v$, that is, if there are words $\alpha$ and
$\beta$ over $A$ and a rule
$(\lambda,\rho)\in R$ such that $u = \alpha \lambda \beta$ and
$v=\alpha \rho\beta$. \textit{Thue equivalence} is the equivalence
relation on $\Astar$ generated by elementary reductions.

There is a multiplication in $\Astar$ given by concatenation. This
induces a multiplication on the set of Thue equivalence
classes. We will work with rules where the set of equivalence classes
is isomorphic to the group $G$.

By no means every set of rules can be used to find the
\shl normal form of a word constructively. We now discuss the
various properties that a set of rules should have in order that reduction
to an irreducible always
gives the \shl normal form of a word.
First we give the assumptions that we will always make about every set
of rules we consider. When constructing a new set of rules, we will
always ensure that these assumptions are correct for the new set.

\subsection{Standard assumptions about rules.}\label{Standard
assumptions about rules.}
\begin{enumerate}
\item\label{reduction to epsilon}[Condition]
For each $x\in A$, $x.\inverse
(x)$ is Thue equivalent to the trivial word $\epsilon$.
The preceding condition
is enough to ensure that the set of Thue equivalence classes is a
group.
If $r=s$ is a defining relation for
$G$, then $r$ is Thue equivalent to $s$.
This ensures that the group of Thue equivalence classes is a quotient
of $G$.
\item\label{same element}[Condition]
If $(u,v)$ is a rule of $R$, then $u>v$ and $\overline u =
\overline v \in G$.
This ensures that the group of Thue equivalence classes is isomorphic
to $G$.
\end{enumerate}

\subsection{Confluence.}\label{Confluence.}[Condition]
This property is one which we certainly desire, but which is hard
to achieve. Given $w$,
there may be different ways to reduce $w$. For example we could look in
$w$ for the first subword that is a left-hand side,
or for the last subword, or just
look for a left-hand side which is some random subword of $w$.
We say that $R$ is \textit{confluent} if the result of fully reducing
$w$ gives an irreducible that is independent of which elementary
reductions were used.

\begin{lemma}\label{properties imply group}[Lemma]
If a set $R$ of rules satisfies the conditions of
\ref{Standard assumptions about rules.} and \ref{Confluence.} then
the set of $R$-irreducibles is mapped bijectively to $G$ and
multiplication corresponds to concatenation followed by reduction.
Under these assumptions,
an $R$-irreducible is in \shl normal form, and conversely;
moreover, each
Thue equivalence class contains a unique irreducible.
\end{lemma}
\begin{proof} The homomorphism $\Astar\to G$ is surjective and, by
\thref{same element}, elementary
reduction does not change the image in $G$. It follows that the
induced map from the set of irreducibles to $G$ is surjective.
Suppose $u$
and $v$ are irreducibles such that $\overline u = \overline v \in G$.
Then $\overline{u .\inverse(v)} = 1_G$. Therefore $u.\inverse(v)$ is
equal in the free group generated by $A$ (with $\inverse(x)$ equated
to the formal inverse of $x$, for each $x\in A$)
to a word $s$ which is
a product of formal conjugates of the defining
relators. Now $u.\inverse(v)$ and $s$ reduce to the same
word, using only reductions that replace $x.\inverse(x)$, where $x\in A$,
by the trivial word $\epsilon$.
By Condition~\ref{reduction to epsilon}, $s$ can be reduced to
$\epsilon$. It follows from Condition~\ref{Confluence.} that
$u.\inverse(v)v$ can be reduced to $v$. It can also be reduced to $u$,
using Condition~\ref{reduction to epsilon} again, and the fact that
$\inverse:A\to A$ is an involution.
It follows from Condition~\ref{Confluence.} that $u=v$, as required.

The description of the multiplication of irreducibles follows from the
fact that multiplication in $\Astar$ is given by concatenation
and the fact that the map $\Astar\to G$ is a homomorphism of monoids.

Since reduction reduces the \shl order of a word, a word in
\shl least normal form must be $R$-irreducible. Conversely, if $u$
is $R$-irreducible, let $v$ be the \shl normal form of $\overline
u$. Then $v$ is also $R$-irreducible, as we have just pointed out,
and $u$ and $v$ represent the same element of $G$. Since the map from
irreducibles to $G$ is injective, we deduce that $u=v$. Therefore $u$
is in \shl normal form.

To show that each Thue equivalence class contains a unique
irreducible, we note that if there is an elementary reduction of $u$
to $v$, then, in case of confluence, any reduction of $u$ gives the
same answer as any reduction of $v$.
\end{proof}

\subsection{Recursive sets of rules.}\label{Recursive sets of
rules.}[Condition]
Another important property (lacked by some of the sets of rules we
discuss)
is the condition that the set of rules be a recursive
set. As opposed to the usual setup when discussing rewrite systems,
we do not require $R$ to be a finite set of
rules---in fact, in this paper $R$ will normally be infinite.
To say that $R$ is \textit{recursive} means that
there exists a Turing machine which can decide whether or
not a given pair $(u,v)$ belongs to $R$.

\begin{definition}\label{U}[Definition]
We denote by $U$ the
set of all rules of the form $(u,v)$, where $u > v$ and $\overline u
= \overline v\in G$. $U$ is called the
\textit{universal set of rules}.
Note that a word is $U$-irreducible if and only if it is in
\shl normal form.
\end{definition}

\begin{lemma}\label{rules and word problem}
The existence of a set of rules $R$ satisfying the conditions of
\ref{Standard assumptions about rules.},
\ref{Confluence.} and \ref{Recursive sets of rules.}
is equivalent to the solubility of the word problem in $G$ and in
this case $U$ defined in \ref{U} is such a set of rules.
\end{lemma}
\begin{proof}
On the one hand,
if we have such a set $R$, then we can solve the word problem by reduction---%
according to Lemma~\ref{properties imply group} a word $w$ reduces to the
trivial word if and only if $\overline w = 1_G$.

On the other hand, if the word problem is solvable, then
the set $U$ of Definition~\ref{U}
is recursive. The various conditions on a set of
rules follow for $U$.
\end{proof}

$U$ can be difficult to manipulate,
even for a very well-behaved
group $G$ and a finite ordered set $A$ of generators, and we therefore
restrict our attention to a much smaller subset, namely the set of
$U$-minimal rules, which we now define.

\begin{definition}\label{minimal}[Definition]
Let $R$ be a set of rules for a group $G$ with generators
$A$. We say that a rule $(u,v)\in R$ is \textit{$R$-minimal} if $v$ is
$R$-irreducible and if every proper subword of $u$ is $R$-irreducible.
\end{definition}

\begin{proposition}\label{conditions and minimals}[Proposition]
\par
\begin{enumerate}
\setlength\itemsep{0pt}
\item\label{conditions minimals satisfy}
The set of $U$-minimal rules satisfies the conditions of
\ref{Standard assumptions about rules.} and \ref{Confluence.}.
In particular they are confluent.
\item\label{length conditions}
Let $(u,v)$ be a $U$-minimal rule and let $u=u_{1}\ldots
u_{n+r}$ and $v=v_{1}\ldots v_{n}$. Then the following must hold:
$0\le r \le 2$; if $n>0$,
$u_{1} \neq v_{1}$; if $n>0$, then $u_{n+r}\neq v_{n}$;
if $r=0$ and $n>0$, then $u_{1} > v_{1}$;
if $r=2$ and $n>0$, then $u_{1}<v_{1}$ and $u_{2} < \inverse(u_{1})$;
if $r=2$ and $n=0$, then $u_1 \le \inverse(u_2)$ and $u_2\le
\inverse(u_1)$.
\item\label{minimals recursive}
The set of $U$-minimal rules
is recursive if and only if $G$ has a solvable word problem.
\end{enumerate}
\end{proposition}
\begin{proof}
If $w$ is $U$-reducible, let $u$ be the shortest prefix of $w$ which
is $U$-reducible. Then every subword of $u$ which does not contain the
last letter is $U$-irreducible.
Let $v$ be the shortest suffix of $u$ which is
$U$-reducible. Then every proper subword of $v$ is $U$-irreducible.
Let $s$ be the \shl normal form for $v$.
Then $(v,s)$ is a $U$-minimal rule. Replacing $v$ in $w$ by $s$ gives
an elementary reduction by a $U$-minimal rule. It follows that reduction
of $w$ using only $U$-minimal rules eventually gives us a $U$-irreducible
word, and this must be the \shl normal form of $w$.
Therefore the conditions of \ref{Standard assumptions about rules.} and
\ref{Confluence.} are satisfied by the set of $U$-minimal rules.

We now prove \ref{length conditions}. Since $u>v$ in the \shl
order, $|u|\ge|v|$. So $r\ge0$. If $r>2$, then $\overline u =\overline
v$ gives rise to $\overline{u_{2}\ldots u_{n+r}} =
\overline{\inverse(u_{1})v_{1}\ldots v_{n}}$. Therefore $u_{2}\ldots
u_{n+r}$ is not in \shl normal form. It follows that
$u_{2}\ldots u_{n+r}$ is $U$-reducible. Therefore $(u,v)$ is
not $U$-minimal. Similar arguments work for the other cases.
This completes the proof of \ref{length conditions}.

Clearly $U$-minimality of a rule can be detected by a Turing machine if the
word problem is solvable. Conversely, if the set of $U$-minimal rules
is recursive, then the word problem can be solved by reduction using
only $U$-minimal rules.
\end{proof}

Now we have a uniqueness result for the set of minimal rules.
\begin{lemma}\label{uniqueness}
Let $R$ satisfy the conditions of \ref{Standard assumptions about rules.}
and \ref{Confluence.}. Suppose every rule of $R$ is
$R$-minimal. Then $R$ is equal to the set of $U$-minimal rules.
\end{lemma}
\begin{proof}
By Lemma~\ref{properties imply group}, the $R$-irreducibles are the same
as the words in \shl normal form.
Let $(u,v)$ be a rule in $R$.
Then $v$ is $R$-irreducible and therefore in \shl normal form.
Also every proper subword of $u$ is in \shl normal form.
Therefore $(u,v)$ is in $U$ and is $U$-minimal.

Conversely, suppose $(u,v)$ is $U$-minimal.
Then $v$ is the \shl normal form of $\overline u$.
By Lemma~\ref{properties imply group} for $R$, $u$ must be $R$-reducible.
Every proper subword of $u$ is already in \shl normal form.
It follows that there is a rule $(u,w)$ in $R$.
Since this rule is $R$-minimal, $w$ is $R$-irreducible.
Therefore $w$ is the \shl normal form of $\overline u$.
It follows that $v=w$. Therefore every $U$-minimal rule is in
$R$.
\end{proof}

We are interested in those pairs $(G,A)$, where $G$ is a group and $A$ is
an ordered set of generators, such that the set of $U$-minimal rules is
not only recursive, but is in fact regular. We now explain what we
mean by \textit{regular} in this context.

We recall that a subset of $\Astar$
is called \textit{regular} if it is equal to $L(M)$, the language accepted by
some finite state automaton over $A$.
(See Definition~\ref{FSA}, where finite state automata are discussed.)
We need to formalize what it means for
an automaton to accept pairs of words over an alphabet $A$.
If the pair of
words is $(abb,ccdc)$, then we have to \textit{pad} the shorter of the two
words to make them the same length, regarding this pair as
the word of length four $(a,c)(b,c)(b,d)(\$,c)$.
In general, given an arbitrary pair of words $(u,v) \in \Astar\times\Astar$,
we regard this instead as a word of pairs by adjoining a {\it padding symbol}
$\$$ to $A$ and then ``padding'' the shorter of $u$ and $v$ so that
both words have the same length.
We obtain a word over $A\cup \{\$\} \times A\cup\{\$\}$.
The alphabet $A\cup\{\$\}$ is denoted $A^+$
and is called the {\it padded extension} of $A$.
The result of padding an arbitrary pair $(u,v)$
is denoted $(u,v)^+$. A word $w \in (A^+)\uast\times(A^+)\uast$
is called {\it padded} if there exists $u,v \in A\uast$
with $w = (u,v)^+$ (that is, at most one of the two components
of $w$ ends with a padding symbol and there are no padding symbols in
the middle of a word).

A set $R$
of pairs of words over $A$ is called {\it regular} if the corresponding
set of padded words is a regular language over the product
alphabet $A^+\times A^+$.
We say that $R$ is accepted by a two-variable finite state automaton over $A$.

\begin{theorem}\label{automatic implies regular}
Let $G$ be a group and let $A$ be a finite set of
generators, closed under taking inverses.
If $(G,A)$ is \shl automatic, then the set of $U$-minimal rules is regular.
\end{theorem}

Having a finite confluent
set of rules does not imply \shl automatic. A counter-example is
given in \cite[page 118]{WPiG}. So the converse of this theorem is not
true.

\begin{proof}
Since we have a \shl automatic structure, the set $L$ of
\shl normal forms is a regular language.
If $x\in A$, the automatic structure includes the multiplier $M_x$,
which is a two-variable automaton over $A$. The language $L(M_x)$ is
the set of pairs $(u,v)$, such that $u,v\in L$ and $\overline{ux}=
\overline v$.
It is not hard to construct from the union of the
$M_x$ an automaton whose language
$P$ is the set of $(u,v)$ such that
$\overline u = \overline v \in G$, $u\in L.A$ and $v\in L$.

We know that
$ (L.A \cap A.L)\cap (\Astar\setminus L)$
is a regular language. Clearly,
this is the set of left-hand sides of $U$-minimal
rules, since it is the set of $U$-reducible
words such that each proper subword is $U$-irreducible.
The set of pairs $(u,v)\in P$, such that $u$ is a left-hand side of a
$U$-minimal rule is easily seen to be the set of all $U$-minimal rules.
\end{proof}

\subsection{Question.}
Suppose $(G,A)$ has a finite confluent set $R$ of \shl reducing
rules which define $G$.
Then it is easy to construct from this a finite confluent set $R'$ of
$R'$-minimal rules defining $G$.
The method is to use minimization, as described in \ref{minimization}.
This set of rules is equal to the set of $U$-minimal rules by
\thref{uniqueness}.

Suppose now that $(G,A)$ has an infinite confluent set $R$ of \shl-reducing
rules defining $G$, and this set is regular. Is the set of
$U$-minimal rules also regular? We know that it is confluent and
recursive by \thref{conditions and minimals}, since $R$ provides a
solution to the word problem.

If $R$ contains all $U$-minimal rules,
then the answer is easily seen to be \textit{yes}. The answer is not
clear to us if $R$ does not contain all minimal rules. There is no loss
of generality in making $R$ smaller so that each proper subword of
each left-hand side is irreducible. But we see no way of
changing $R$ so as to ensure that each right-hand side is irreducible,
while maintaining $R$'s property of being regular.

\subsection{Objective.}\label{Objective.}
In this paper we present a procedure which,
given a set of rules satisfying
the conditions of
\ref{Standard assumptions about rules.}, changes the set of rules so that it
becomes ``more confluent''. More precisely, the set of words for
which all reductions give the same irreducible, and this irreducible
is in \shl normal form, increases with time. If we fix attention on
a single word this will eventually be included in the set. However, in
general, because of the insolubility of the word problem, it is not
in general possible to know when that time has been arrived at.

For a group where the set of all $U$-minimal rules
(see Definition~\ref{U}) is the set of all pairs accepted by a
two-variable minimal PDFA $M$ (these concepts are defined in
\ref{FSA}), our procedure gives
rise to $M$ after a finite number of steps.

For many undecidable problems, there is a ``one-sided'' solution. The
technical language is that a certain set is recursively enumerable,
but not recursive. For
example, consider a fixed group for which the word problem is undecidable.
Given a word $w$ in the generators, if you are correctly
informed that $\overline w = 1_{G}$, then this can be verified by a
Turing machine. All that you have to do is to enumerate products of
conjugates of the defining relators, reduce them in the free group
on the generators, and see if you get $w$, also reduced in the free
group. If $w$ represents the identity then you will prove this
sooner or later. If it's not the identity, the process continues for
ever.

We know that there is no algorithm which has as input a finite
presentation of a group and outputs whether the group is trivial or
not (see \cite{Rabin:unsolvability}).
It follows easily that there is no algorithm which has as input a
finite presentation and outputs either an FSA accepting the set of
$U$-minimal rules or correctly answers \textit{There is no such FSA}.
For, in the case of the trivial group,
the set of $U$-minimal rules is finite---for each element $x\in A$, we
have the rule $(x,\epsilon)$---and so it is certainly regular.

But the situation is even worse than this. We do not even know of a
one-sided solution to the problem of whether the set of $U$-minimal
rules is regular. If the set of $U$-minimal rules is regular,
our procedure will eventually produce a candidate
with some indication that it is correct, but we will not know \textit{for sure}
whether the answer is correct or incorrect.

What is at issue is whether there is an algorithm which has as its
input a regular set of \shl rules for a group
and outputs whether or not the set of rules is confluent.
For finite sets of rules the question of confluence is
decidable by classical critical pair analysis
which we describe in \thref{Standard Knuth--Bendix}.
However, for
{\it infinite} rewriting systems the confluence question is, in
general,
undecidable. Examples exhibiting undecidability are given in
\cite{InfiniteRegularThue}.
They are length-reducing rewriting
systems $R$ which are regular in a very strong sense:
$R$ contains only a finite
number of right-hand sides and for each right-hand side $r$, the set
$\{l : (l,r) \in R\}$ is a regular language.
These examples are in the context of rewriting for monoids.
As far as we know,
there is no known example of undecidability if we add to the hypothesis
that the monoid defined by $R$ is in fact a group.

In the special case where $(G,A)$ is \shl automatic, there
is a test for confluence of a set of rules satisfying the conditions
of \ref{Standard assumptions about rules.},
namely the axiom-checking procedure described in theory in
\cite{WPiG} and carried out in practice in Derek Holt's
\textit{kbmag} programs \cite{Holt:KBMAG}.
%
%
%
%
%
%
%
%
%
%
%
%

%% file: Sec3.tex
\section{Welding}
\label{Welding}[Section]

In this section we start with an example which motivates the operation
of welding. We then give a formal definition, and prove that the
operation gives rise to a function from the set of regular languages
to the set of regular languages. We then define the concept of a rule
automaton---this is a finite state automaton in two variables which
can recognize when certain words in the generators are equal in the
associated group.
We show that a welded rule automaton is also a rule automaton.

\subsection{A motivating example.}\label{motivating example}
We will use the standard generators $x$, $y$,
and their inverses $X$ and $Y$ for the free abelian group on two
generators. We will impose different orderings on this set of four generators,
and, as described in \ref{Objective.}, see what kind of confluent sets 
of rules emerge.

Consider the alphabet $A=\{x,X,y,Y\}$ with the ordering $ x < X < y < Y$,
and denote the identity of $A\uast$ by $\epsilon$.
Let $R$ be the rewriting system on $A\uast$
defined by the set of rules
$$\{(xX,\epsilon),\ (Xx,\epsilon),\ (yY,\epsilon),\ (Yy,\epsilon),
\ (yx,xy),\ (yX,Xy),\ (Yx,xY),\ (YX,XY)\}.$$
It is straightforward to see
that $R$ is a confluent system.

We now change the ordering of the set of generators to
$x < y < X < Y$ and correspondingly interchange the sides of the sixth
rule getting $(Xy,yX)$ and an order reducing set of rules.
Once again the rules define the free abelian group on two generators.
But this time there can be no finite confluent set of rules.
To see this, we consider the set of words
$\{xy^nX : n \in \naturals\}$. None of these is in \shl normal form.
By \thref{properties imply group}, each of these words is reducible relative
to any confluent set of rules.
On the other hand,
each proper subword of one of
the words $xy^nX$ is clearly in \shl normal form and is therefore
irreducible. It follows that a confluent
set of rules must contain each of the words $xy^nX$
as a left-hand side.
In this situation, the classical Knuth--Bendix procedure (see
\thref{Standard Knuth--Bendix})
will never terminate,
and the same is true for any method of which generates only a finite 
number of rules at each step.

We will now introduce a new procedure, which we call
\textit{welding}. This can produce an infinite set of rules
from a finite set of rules in a finite number of steps.
Welding is central to the main procedure of the computer
program described in this paper.

First we need to give some standard definitions.
\begin{definition}\label{FSA}[Definition]
A finite state automaton (abbreviated FSA) $M$ over a finite alphabet
$A$ is a finite graph
with directed edges and the following additional properties.
Each edge (called an \textit{arrow} in this
context) is either labelled with an element of
$A$ or is unlabelled. Unlabelled arrows are sometimes labelled with
$\epsilon$, which stands for the empty word, and are called
\textit{$\epsilon$-transitions}.
The vertices of the graph are called \textit{states}.
Some of the states are labelled as \textit{initial states} and some as
\textit{final states}. The language $L(M)$ accepted by $M$ is the set
of words over $A$ which are traced out by paths of arrows
which start at some initial state and end at some final state.
An FSA is said to be \textit{partially deterministic} (abbreviated 
PDFA) if it has no
$\epsilon$-transitions, if there is exactly one initial state and if,
for each state $s$ and each $x\in A$, there is at most one arrow from $s$
with label $x$. An FSA is said to be \textit{trim} if, for each state
$s$, there is a path of arrows which starts at an initial state,
and ends at a final state, with $s$ lying on the path.
The \textit{reversal} of a finite state automaton is the same graph with the
same labelling, but with each arrow reversed, with each initial state
changed to be a final state and each final state changed to be an
initial state.
A \textit{non-deterministic} automaton NFA is an automaton with
$\epsilon$-transitions and/or some states $s$ having more than one
arrow from $s$ having the same label.
\end{definition}

\begin{definition}
\label{welding}
An FSA is called {\it welded} if it is partially deterministic, trim and has a
partially deterministic reversal.
These conditions imply that, given $x\in A$ and a state $t$, there is at
most one $x$-arrow with target $t$ and also that there is exactly one
initial state and one final state.
\end{definition}

Given a trim non-empty FSA $M$, we can form a welded automaton from it as 
follows.
Given any $\epsilon$-arrow $(s,\epsilon,t)$, we may identify $s$ with
$t$.
Given distinct initial states $s_1$ and $s_2$, we may identify $s_1$
with $s_2$.
Given distinct final states $t_1$ and $t_2$, we may identify $t_1$
with $t_2$.
Given distinct arrows $(s,x,t_1)$ and $(s,x,t_2)$, we may identify
$t_1$ with $t_2$.
Given distinct arrows $(s_1,x,t)$
and $(s_2,x,t)$, we may identify $s_1$ with $s_2$.
Immediately after any identification of two states, we change the set of
arrows accordingly, omitting any $\epsilon$-arrow from a state to itself.
Since the number of states continually decreases,
this process must come to an end, and at this point the automaton is
welded.

\subsection{Welding in our example.}
Let us see how this works on the example given in \ref{motivating 
example}. For the moment we won't try to justify the correctness of 
our procedure, that is, that the new rules that welding produces are 
valid rules; we will just carry out the procedure to show how it works.
Justification comes from the consideration of rule automata---see \thref{validrules}.

\begin{figure}[htbp]
{\fontsize{10}{12pt}\selectfont
\begin{center}\begin{picture}(240,50)(20,0)
\put(-12.5,20){\vector(1,0){10}}
\put(0,20){\circle{5}}
\put(-2,5){$1$}
\put(2.5,20){\vector(1,0){40}}
\put(10,30){$(x,y)$}
\put(45,20){\circle{5}}
\put(42.5,5){$2$}
\put(47.5,20){\vector(1,0){40}}
\put(55,30){$(y,y)$}
\put(90,20){\circle{5}}
\put(87.5,5){$3$}
\put(110,20){$\ldots$}
\put(140,20){\circle{5}}
\put(137.5,5){$n$}
\put(142.5,20){\vector(1,0){40}}
\put(150,30){$(y,y)$}
\put(185,20){\circle{5}}
\put(173.5,5){$n+1$}
\put(187.5,20){\vector(1,0){40}}
\put(195,30){$(y,\$)$}
\put(230,20){\circle{5}}
\put(218.5,5){$n+2$}
\put(232.5,20){\vector(1,0){40}}
\put(275,20){\circle*{5}}
\put(263.5,5){$n+3$}
\put(240,30){$(X,\$)$}
\end{picture}
\end{center}
}
\caption{\sf The PDFA $M(r_n)$ for $n>2$.}
\label{ruletopdfa}
\end{figure}
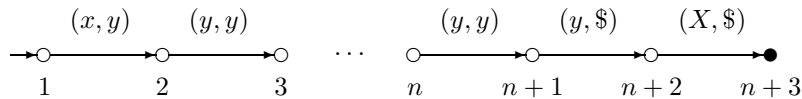

We consider the rule $r_n = (xy^nX,y^n)$
for some $n\in \naturals$. The corresponding padded word $r_n^+$ gives
rise to an $(n+3)$-state PDFA $M(r_n)$ whose accepted language consists solely
of the rule $r_n$. For $n>2$ this PDFA is shown in Figure~\ref{ruletopdfa}.

\begin{figure}[htbp]
{\fontsize{10}{12pt}\selectfont
\begin{center}
\begin{picture}(220,60)
\put(0,17.5){\circle{5}}
\put(-2.5,2.5){$1$}
\put(-12.7,17.5){\vector(1,0){10}}
\put(2.4,17.5){\vector(1,0){70}}
\put(74.8,17.5){\circle{5}}
\put(72.3,2.3){$2$}
\put(74.8,30){\circle{20}}
\put(71,37){$<$}
\put(64,46){$(y,y)$}
\put(22,22){$(x,y)$}
\put(77.2,17.5){\vector(1,0){70}}
\put(149.6,17.5){\circle{5}}
\put(147.1,2.5){$3$}
\put(105,22){$(y,\$)$}
\put(152.1,17.5){\vector(1,0){70}}
\put(224.6,17.5){\circle*{5}}
\put(221.1,2.5){$4$}
\put(179,22){$(X,\$)$}
\end{picture}
\end{center}
}
\caption{\sf A PDFA isomorphic to $\Weld{M_n}, n > 1$.}
\label{weldrules}
\end{figure}
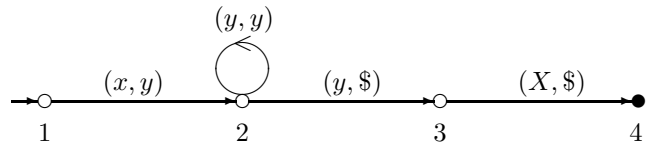

Continuing the discussion of the rules for a free abelian group on two
generators, we define $M_n$ to be the disjoint union
$\bigcup \{M(r_1),\ldots,M(r_n)\}$ of the
automata $M(r_1),\ldots,M(r_n)$, with set of initial (final) states
equal to the collection of initial (final) states for the various $M(r_i)$.
If $n>1$ then $\Weld{M_n}$ is
isomorphic to the PDFA given in Figure~\ref{weldrules}, and the accepted language
of this PDFA is the set of rules $\{r_i : i \in \naturals\}$. This is
independent of $n$ if $n>1$.

So in this
example, after only two steps,
the welding procedure provides us with a PDFA whose accepted
language consists of an infinite set of identities between words
in the free abelian group.
Moreover, by using this PDFA to
define a suitable reduction procedure,
each of the words $xy^nX$ with $ n\in \naturals$
can be reduced to the \shl normal form.

For this group with the given ordering on the generators, it is not
hard to show that by welding the original defining rules for the group together
with the $4$ rules $\{(xyX,y),(xy^2X,y^2),(yXY,X),(yX^2Y,X^2)\}$, we obtain a
PDFA whose accepted language is a confluent set of rules (provided we
adjust the automaton to ensure that only padded pairs of words $(u,v)^+$ are
accepted, with $u>v$).
Any reduction procedure using this infinite set of rules
will reduce {\it any} word to its \shl normal form.

The next theorem is a general result about the welding of finite state
automata which need have nothing to do with groups. It's a result
which is reassuring, but, logically, it is entirely unnecessary for
understanding other parts of this paper. Readers pressed for time
should skip it.
\begin{theorem}\label{welding unique}
Given a trim non-empty FSA $M$, all welded automata obtained from it
as above $($no matter in what order the states and arrows are identified
to each other$)$ are the same, except that the names of the states may
be different.
The automaton $Q$ thus obtained is a minimal PDFA and $Q$
depends only on the language $L(M)$, up to changing the names of the
states.
It follows that
welding can be regarded as an operation on regular languages,
independent of the automaton used to encode them.
\end{theorem}

\begin{proof}
For each $x\in A$, let $x^{-1}$ be its formal inverse and let $A^{-1}$
be the set of these formal inverses. We form from
$M$ an automaton over $A\cup A^{-1}$ by adjoining an arrow of the form
$(t,x^{-1},s)$ for each arrow $(s,x,t)$ of $M$, and
adjoining an arrow
$(t,\epsilon,s)$ for each arrow $(s,\epsilon,t)$ unless it's already
there. We also adjoin
$(s_1,\epsilon,s_2)$ if $s_1$ and $s_2$ are either both initial states
or both final states, unless these arrows are already there.
We denote this new automaton by $N$.
$N$ has the same initial and final states as $M$.

Let $F$ be the free group generated by $A$.
We define a relation on the set of states of $N$ by $s \sim t$ if
there is a path of arrows from $s$ to $t$ in $N$ whose label gives the
identity element of $F$. This is clearly an equivalence relation.
Let $Q$ be the automaton defined as follows. Each state of $Q$ is one of the
equivalence classes above.
The unique initial state of $Q$ is the unique equivalence class containing
all initial states of $N$. The unique final state of $Q$ is the unique
equivalence class containing all final states of $N$.
Let $S$ be one equivalence class and
$T$ another, and let $x\in A$.
We have an arrow $x:S\to T$ in $Q$ if there is an $s\in S$ and a $t\in T$
and an arrow $x:s\to t$ in $M$.
It is easy to see that $Q$ is welded, and it follows that it is a
partial deterministic automaton.

If $M$ starts out by being welded, then it is easy to see that $Q=M$,
up to the naming of states.

Consider the identifications of states and arrows
made during welding (see the passage following \thref{welding}).
Let $M=M_0,M_1,\ldots,M_k$ be the sequence of automata obtained by
identifying at each step only one state with another state or deleting
one arrow labelled $x$ from a state $s$ to state $t$ if there are several
arrows labelled $x$
from $s$ to $t$ or deleting one $\epsilon$-arrow from a state to itself.
Here $M_k$, the last automaton in the list, is a welded automaton.

We assign to each state $s$ of $M_i$
the set of all states of the original automaton $M$ which are
identified to make $s$.
A state $q$ of $Q(M_i)$ is a set of states of $M_i$, and this is a set
of subsets of the state set of $M$.
By taking the union, we can instead regard $q$ as a set of states of
$M$.
This loses some of the structure, but only an irrelevant part.

With this interpretation, we see that the states of $Q(M_i)$ are
identical to those of $Q(M_{i+1})$.
Moreover, all arrows in $Q(M_i)$ are inherited from $M$ via $M_i$.
It follows that the automaton $Q(M_i)$ is independent of $i$.
So we have $Q= Q(M)= Q(M_k) = M_k$.
This shows that $Q$ is independent of the order in which the
identifications are carried out. In fact $Q$ can be characterized as
the largest welded quotient of $M$.

We claim that every element of $L(Q)$ arises as follows, and that only
elements of $L(Q)$ arise in this way. Let $(w_1,
w_2, \ldots, w_{2k+1})$ be a $(2k+1)$-tuple of elements of $L(M)$, where
$k\ge 0$. Now consider
$$w_1 w_2^{-1} \ldots w_{2k}^{-1} w_{2k+1} \in F,$$
and write it in reduced form, that is, cancel adjacent formal
inverse letters wherever possible. If the result is in $A\uast$, that
is, if after cancellation there are no inverse symbols,
then it is in $L(Q)$.

To prove this claim, we proceed as follows.
For each state $s$ of $M$,
we fix a path of arrows $p_s$ in $M$ from an initial state to $s$ and
a path of arrows $q_s$ from $s$ to a final state.
If $s$ is an initial state, we define $p_s$ to be the trivial path.
If $s$ is a final state, we define $q_s$ to be the trivial path.

Start with an arbitrary element $w\in L(Q)$.
We must show that $w$ can be produced in the way described above.
Now $w$ is the label of a path of arrows in $Q$, starting from the
initial state of $Q$ and ending at the final state of $Q$.
Recalling the definition of a state of $Q$, we
can replace this path by a path of arrows in $N$,
which alternately traverses
a path of arrows in $N$ labelled by a word over $A\cup A^{-1} \cup\{\epsilon\}$
which reduces to the identity element in $F$,
and an arrow of $N$ labelled by a letter in $w$.
The path in $N$ starts at an initial state of $N$ and ends at a final state
of $N$.
We write the path as a composite of arrows $u_i$ in $N$.

If $u_i:s\to t$ is an arrow in $M$,
we replace it by $p_s^{-1}\left(p_s u_i q_t\right)q_t^{-1}$.
Otherwise, if the inverse of $u_i:s\to t$ is an arrow of $M$, we
replace $u_i$ by $q_s \left(q_s^{-1} u_i p_t^{-1}\right) p_t$.
(We consider the inverse of an $\epsilon$-arrow to be an
$\epsilon$-arrow.)
Otherwise $s$ and $t$ are both initial states or both final states and
$u_i$ is an $\epsilon$-arrow and we leave $u_i$ unaltered.

Each expression within parentheses in the preceding paragraph
therefore give either some $w_i\in
L(M)$ (possibly empty) or the formal inverse of such a word.
Outside these parentheses we obtain expressions like
$\epsilon$, $q_s^{-1}q_s$, $p_s p_s^{-1}$, $p_sq_s$ or $q_s^{-1}p_s^{-1}$.
In the first three cases, we omit the expressions. In the last two
cases, the expression represents either  $w_i\in L(M)$,
or the formal inverse of such a word.
The path starts at an initial state of $N$ and ends at a final state.
So, if the set of initial states is disjoint from the set of final
states, then the expression of $w$ as a product in the free group $F$
of elements of $L(M)$
and their formal inverses must have an odd number of factors.
If the set of initial states meets the set of final states, then the
trivial word is an element of $L(M)$, and we can use this to make sure
that the number of factors is odd. This completes the claim in one
direction.

Conversely, suppose we are given the $w_i \in L(M)$ as in the claim.
Then $w_i$ is the label on a path of arrows in $M$ from an initial
state to a final state.
By inserting $\epsilon$-arrows in $N$ to join initial states or to
join final states, we find that $w_1 w_2^{-1} \ldots w_{2k}^{-1} w_{2k+1}$
is the label of a path of arrows in $N$ from an initial state to a
final state. An elementary cancellation in $F$ corresponds to the fact that
two states of $N$ give rise to the same state of $Q$.
Carrying out all the elementary cancellations possible, if we are left
only with a word over $A$, we have defined a path of arrows in $Q$
from the initial state of $Q$ to the final state of $Q$.
So we have found an element of $L(Q)$, as claimed.

A welded automaton is minimal.
For let $s$ and $t$ be distinct states, and let $u$ and $v$ be words
over $A$ which lead from $s$ and $t$ respectively to the unique final
state. Then $u$ does not lead from $t$ to the final state and $v$ does
not lead from $s$ to the final state (otherwise $s$ and $t$ would be
equal). It follows that $s$ and $t$ remain distinct in the minimized
automaton.
\end{proof}

If $M$ is a non-empty trim FSA, we denote by $\Weld{M}$ the PDFA
obtained from it by welding.
To compute $\Weld{M}$ efficiently, we first
add ``backward arrows'' to $M$.
That is, for each arrow $(s,x,t)$ in $M$, including $\epsilon$-arrows,
we add the arrow $(t,x',s)$, where $x'$ represents a backwards version
of $x$. We also add $\epsilon$-arrows to connect the initial states,
and $\epsilon$-arrows to connect the final states.
We then make use of a slightly modified version of
the coincidence procedure of Sims given in
\cite[4.6]{sims}. When this stops we have a welded automaton.

In practice, in the automata which we want to weld, backward arrows are
needed in any case for some algorithms which we need.
The procedure described
in the preceding paragraph therefore fits our needs particularly well.

For the welding procedure to be used in a general Knuth--Bendix situation, we
need to show that any rules obtained are valid identities in the corresponding
monoid. We now show that if the monoid is a group (the situation
we are interested in), any rules obtained are valid identities.

\begin{definition}
\label{worddiff}[Definition]
Let $A$ be a finite inverse closed set of monoid generators for a group
$G$ and, as before,  denote images under the surjection $(A^+)\uast
\rightarrow G$ by overscores. A {\it rule automaton for $G$} is a
two-variable FSA $M=(S,A^+~\times~A^+,\mu,F,S_0)$ together with a
function $\phi_M : S \rightarrow G$ satisfying
\begin{enumerate}
\item $F,S_0 \neq \emptyset$.
\item If $s$ is an initial or final state then $\phi_M(s)=1_G$.
\item \label{worddiff:3}
For any $s,t \in S$ and $(x,y) \in A^+ \times A^+$ with
$(s,(x,y),t)\in\mu$
we have $\phi_M(t)~=~\overline{x}^{-1}\phi_M(s)\overline{y}$.
\item For any $s,t\in S$ with $(s,\epsilon,t)\in\mu$ we have
$\phi_M(s)=\phi_M(t)$.
\end{enumerate}
\end{definition}


\begin{example}
\label{ruleworddiff}
If $A$ is a finite inverse closed set of monoid generators for a group $G$ and
$r = (u,v)\in A\uast \times A\uast$ satisfies $\overline u = \overline v$ then,
as in Figure~\ref{ruletopdfa}, writing $r^+$ as a word
$(u_1,v_1)\cdots(u_n,v_n)\in(A^+\times A^+)\uast$, we obtain an $(n+1)$-state
rule automaton
$M(r)=(\{s_0,\ldots,s_{n}\},A^+~\times~A^+,\mu,\{s_0\},\{s_n\})$ for $G$ where
the arrows are given by
$$\mu(s_i,(u_{i+1},v_{i+1})) = s_{i+1},\, 0\leq i\leq n-1.$$
The function $\phi=\phi_{M(r)}$ assigning group elements to states is
defined inductively by $\phi(s_0) = 1_G$ and
$\phi(s_i) = \overline {u_i}^{-1} \phi(s_{i-1}) \overline{v_i}$
for $1\leq i \leq n$. As usual, the padding symbol is sent to $1_G$. The
fact that $\overline u = \overline v$ ensures that Condition $2$ of
\thref{worddiff} is satisfied.
\end{example}

\begin{remark}
\label{worddiffpdfa}
For a two-variable FSA $M$ which is a rule automaton,
the PDFA $P$ obtained by applying
the subset construction to the (non-empty) set of initial states of $M$ (and
the sets that arise), is also a rule automaton for $G$, where the map $\phi_P$
is induced from $\phi_M$.
The fact that this map is well-defined follows from Conditions $2,3$ and $4$
of \thref{worddiff} and the fact that $P$ is connected (by construction).

The same remark applies to the modified subset construction described
in Section~\ref{A modified determinization algorithm}.
\end{remark}

\begin{proposition}
\label{validrules}
Let $A$ be a finite inverse closed set of monoid generators for a group $G$ and
suppose that $M$ is a rule automaton for $G$. Then
\begin{enumerate}
\item Every pair $(u,v)\in L(M)$ gives a valid identity $\overline u =
\overline v$ in $G$.
\label{validrules1}
\item $\Weld M$ is a rule automaton for $G$.
\label{validrules2}
\end{enumerate}

Consequently every
accepted rule (that is, an accepted pair $(u,v)$ such that $u>v$)
of $\Weld M$ is a valid identity in $G$.
\end{proposition}

\begin{proof}
To prove \ref{validrules1}, let $r = (u,v)\in A\uast\times A\uast$
be an accepted rule of $M$
and write the padded
word $(u,v)^+$ as $(u_1,v_1)\cdots(u_n,v_n)$.
Then in the PDFA $P$ obtained from $M$ (as in \thref{worddiffpdfa}), there
exists a sequence of states $s_0,\ldots,s_n$ of $P$, such that $s_0$ is the
initial state, $s_n$ a final state, and, for each
$i, 1\leq i\leq n$, there is a arrow from $s_{i-1}$ to $s_i$ labelled by
$(u_i,v_i)$.
Hence, from Condition $3$ of \thref{worddiff}, we have
$$\phi_P(s_i) = {\overline{u_i}}^{-1}\cdots{\overline{u_1}}^{-1}
{\overline{v_1}}\cdots{\overline{v_i}}, \mbox{ for all } i \mbox{ with }
0\leq i\leq n.$$
Condition 2 of \thref{worddiff} tells us that $\phi_P(s_n) = e$.
It follows that $\overline{u_1\cdots u_n} = \overline{v_1\cdots v_n}$,
and therefore the rule $r$ is valid in $G$.

To prove $2$, we need only show that when any of the operations described
just after \thref{welding} is applied to a rule automaton $M$, we
continue to have a rule automaton.
This is obvious. The final statement is now immediate.
\end{proof}

\begin{corollary}
\label{validrulescorollary}
Let $A$ be a finite inverse closed set of monoid generators for a group $G$ and
suppose that $r_1,\ldots,r_m\in A\uast \times A\uast$ give valid identities in 
$G$. Then any rule accepted by $\Weld{M(r_1),\ldots,M(r_m)}$ also
gives a valid identity in $G$.
\end{corollary}

\begin{proof}
For $1\leq k\leq m$ let $M(r_k)$ be the rule automaton for $G$ as in
\thref{ruleworddiff}. Then the disjoint union
$\bigcup \{M(r_1),\ldots,M(r_m)\}$
is also a rule automaton for $G$ and so the result follows by
\ref{validrules}.
\end{proof}

\begin{remark}
\label{identify}
Given a rule automaton $M$ for a group $G$, the map $\phi_M$ may not be
injective. In order to think of the matter constructively, we specify
the values of $\phi_M$ by representing them as words in the
generators.
The undecidability of the word problem implies that the injectivity of
$\phi_M$  might 
be impossible to decide, though sometimes we are in a position to know
whether $\phi_M$ is injective or not.
Even if $\phi_M$ is not injective, the rule
automaton $M$ can still be useful for finding equalities in the group
$G$. $M$ may not tell the whole truth, but it does tell nothing but
the truth.
However, if $\phi_M(s)=\phi_M(t)$ and we can somehow
determine that this is the case, then we can connect
$s$ to $t$ by an $\epsilon$-arrow, and we still have a rule automaton.
If we then weld, $s$ and $t$ will be identified. In this way, with
sufficient investigation, we can hope to make $\phi_M$ injective in
particular cases, even though we know that in general this is an
impossible task.
\end{remark}

\begin{theorem}
Let $G$ be a group and let $A$ be a finite set of generators, closed
under taking inverses. If $G$ is determined by a regular set of
\shl-reducing rules, then $G$ is finitely presented.
\end{theorem}
\begin{proof}
Let $M$ be the finite state automaton accepting the rules in our
regular set.
Then $M$ can be given the structure of a rule automaton, associating
to each state of $M$ a word over $A$.
By \thref{worddiff}, each
arrow $(x,y):s\to t$ in $M$ gives rise to a relation of the form
$\phi_M(t)~=~\overline{x}^{-1}\phi_M(s)\overline{y}$.
There are only a finite number of these, and they can
clearly be combined to prove that $\bar u = \bar v$ for any $(u,v)$
accepted by $M$. It follows that this finite set of relators is a
defining set for $G$.
\end{proof}
%
%
%
%
%
%
%
%
%
%
%
%
%

%% file: Sec4.tex
\section{Standard Knuth--Bendix.}
\label{Standard Knuth--Bendix}[Section]

We recall the classical Knuth--Bendix procedure. Later we will explain
how our procedure differs from it.
We continue to restrict to the \shl case and to groups.
Suppose $G$ is a group given by a finite set of generators and
relators. We define $A$ to be the set of generators together with
their formal inverses. Our initial set of rules
consists of all rules of the
form $(x.\inverse(x),\epsilon)$ for $x\in A$, together with 
all rules of the form $(r,\epsilon)$, where $r$ varies over the
finite set of defining relators for $G$.

After running the Knuth--Bendix procedure (which we are about to
describe) for some time, we will still have a finite set $R$ of rules.
As always, we assume that $R$ satisfies
Conditions~\ref{Standard assumptions about rules.}.

To test for confluence of a finite set of rules, we need only do critical pair
analysis, as explained in \ref{Critical pair analysis.}, \ref{first
case} and \ref{second case}. The proof of this is as follows.

Suppose $R$ is not confluent.
Let $w$ be the \shl least word over $A$ for which there are two
different chains of elementary reductions giving rise to distinct
irreducibles.
Since $w$ is shortest, it is easy to see that the first elementary
reductions in the two chains must overlap.

\subsection{Critical pair analysis.}
\label{Critical pair analysis.}
A pair of rules
$(\lambda_1,\rho_1)$ and $(\lambda_2,\rho_2)$ can overlap
in two possible ways.
First, a non-empty word $z$ may be a suffix of $\lambda_1=s_1z$
and a prefix of
$\lambda_2=zs_2$ (or vice versa).
Second, $\lambda_2$ may be a subword of
$\lambda_1$ (or vice versa) and we write $\lambda_1=s_1\lambda_2s_2$.

These cases are not disjoint. In particular, if one of $s_1$
and $s_2$ is trivial in the second case, it can equally well be
treated under the first case with $z$ equal either to $\lambda_1$ or
to $\lambda_2$.

\subsection{First case of critical pair analysis.}
\label{first case}
In the first case, there are two elementary reductions of
$u=s_1zs_2$, namely to $\rho_1s_2$ and to $s_1\rho_2$.
Further reduction to irreducibles either gives the same irreducible
for each of the two computations, or else
gives us distinct irreducibles $v$ and $w$.
From Conditions~\ref{Standard assumptions about rules.} we deduce that
$v$ and $w$ represent the same element of $G$.
So, if $v$ and $w$ are distinct, we augment $R$ with the
rule $(v,w)$ if $w < v$ or
with $(w,v)$ if $v < w$.
Clearly Conditions~\ref{Standard assumptions about rules.} are
maintained.

Note that it is important to allow $(\lambda_1,\rho_1) = (\lambda_2,\rho_2)$
in the case just discussed,
provided there is a $z$ which is both a proper suffix and
a proper prefix of $\lambda_1=\lambda_2$.

\subsection{Second case of critical pair analysis.}
\label{second case}
In the second case, there are two elementary reductions of
$u=\lambda_1=s_1\lambda_2 s_2$, namely to $\rho_1$ and to $s_1\rho_2s_2$.
If $\rho_1$ and $s_1\rho_2s_2$ reduce to distinct irreducibles
$v$ and $w$, we augment $R$
with either $(v,w)$ or with $(w,v)$, depending on whether $v>w$ or
$w>v$.

\subsection{Omitting rules.}
\label{omitting rules}
In practice, it is important to remove rules which are redundant, as
well as to add rules which are essential. Omitting rules is
unnecessary in theory, provided that we have unlimited time and space at our
disposal. In practice, if we don't omit rules, we are liable to be
overwhelmed by unnecessary computation. Moreover, nearly all programs
in computational group theory suffer from excessive demands for space.
Indeed this is one of the reasons for developing the algorithms and
programs discussed in this paper. So it is important to throw away
information that is not needed and doesn't help.

For this reason, in Knuth--Bendix programs one looks from
time to time at each rule $(\lambda,\rho)$ to see if it can be
omitted. If a proper subword of
the left-hand side can be reduced, then we
are in the situation of \ref{second case}. If the two reductions
mentioned in \ref{second case}
lead to the same irreducible, we omit $(\lambda,\rho)$ from the
set of rules. If the two reductions lead to different irreducibles,
then we augment the set of rules as described in \ref{second case} and
again omit $(\lambda,\rho)$.
We also investigate whether the right-hand
side $\rho$ of a rule $(\lambda,\rho)$ is reducible to $\rho'$.
If so, we can omit $(\lambda,\rho)$ from $R$ and replace
it with the rule $(\lambda,\rho')$.

It is easy to see that such omissions do not change the Thue
equivalence classes.
The process of analyzing critical pairs and augmenting or diminishing
the rule set while maintaining the conditions of \ref{Standard
assumptions about rules.} is called the {\it Knuth--Bendix Process}.

If the Knuth--Bendix process terminates, every left-hand side having
been checked against every left-hand side in critical pair analysis
without any new rule being added, we know that we have a finite
confluent system of rules.
Usually it
does not terminate and it produces new rules {\it ad infinitum}.

\begin{definition}\label{fair}[Definition]
It is important that the process be \textit{fair}.
By this we mean
that if you fix your attention on two rules at any one time, then
either their left-hand sides must have already been, or must
eventually be, checked for overlaps; or one
or both of them must eventually be omitted.
If the process is not fair, it might concentrate exclusively on one
part of the group: for example, in the case of the product of two
groups, the process might pay attention only to one of the factors.
\end{definition}

\subsection{The limit of the process.}\label{limit}
As the Knuth--Bendix process proceeds, $R$ changes and the set of
$R$-reducibles steadily increases. This is obvious when we add a rule
as in \ref{first case} and \ref{second case}. It is also easy to see
when we omit a rule---we need only check that if we omit
$(\lambda,\rho)$ from $R$ as in \ref{omitting rules},
then $\lambda$ remains reducible.

Now let us fix a positive integer $n$.
Eventually the set of reducibles of length at most $n$ stops increasing
with time, and the set of irreducibles of length at most $n$ stops
decreasing.
Since the word problem is in general insoluble, we will in general not know for
sure at any one time or for any fixed $n$ whether the set of
reducibles has stopped increasing. It may look as though it has
permanently stabilized and then suddenly start increasing again. 

Once stabilized, we know by \thref{fair}
that any two reductions of a given word of
length at most $n$ will give the same irreducible (otherwise a new
rule would be added at some time, creating one of more new reducibles
of length at most $n$).
It follows that if
we take the limit of the set of rules (the set of rules which appear
at some time and are never subsequently omitted),
then we have a confluent set of rules. We deduce from
\thref{properties imply group} that, after stabilization of the set of
reducibles of length at most $n$,
any irreducible of length at most $n$ is in \shl normal form.
In fact, at this point, the set of rules with left-hand side of length at most
$n$ coincides with the set of $U$-minimal rules in $U$ (defined in
\ref{U} and \ref{minimal}).

\subsection{Knuth--Bendix pass.}
\label{Knuth--Bendix pass}
One procedure for carrying out the Knuth--Bendix process is to divide the
finite set \Store of rules found so far into three disjoint subsets.
The first subset,
called {\Considered}, is the set of rules whose left-hand sides have been
compared with each other and with themselves for overlaps.
The second set of rules, called {\This}, is the set of rules waiting to
be compared with those in {\Considered}. The third set, called {\New},
consists of those rules most recently found.
Here we only sketch the process. Fuller details of our more elaborate
form of Knuth--Bendix are provided in \thref{our version}.

The Knuth--Bendix process proceeds in phases, each of which is
called a \textit{Knuth--Bendix pass}. Each pass starts by looking at each
rule in {\Considered} and seeing whether it can be deleted as in
\ref{omitting rules}.
Consideration of an existing rule in {\Considered} can lead to a
new rule, in which case the new rule is added to {\New}.

Next, we look at each rule $r$ in \New to see if it is can be omitted
or replaced by a better rule, a process which we call
\textit{minimization}.
The details of our minimization procedure
will be given in \ref{minimization}.
If the minimization procedure changes a rule, the old rule is either
deleted or marked for future deletion. The new rule is added to \This.
Eventually \New is emptied.

We then look at each rule in \This.
Its left-hand side is compared with
itself and with all the left-hand sides of rules in
{\Considered}, looking for overlaps as in \ref{first case}.
Any new rules found are added to {\New}.
Then $r$ is moved into {\Considered}.
Eventually \This becomes empty.

We then proceed to the next pass.

%% file: Sec5.tex
\section{Our version of Knuth--Bendix.}
\label{our version}[Section]

In this section we consider a
rewriting system which is the accepted language of a rule automaton for some
finitely presented group. We call the automaton \Rules.
We  describe a Knuth--Bendix type algorithm for such
a system. In light of the undecidability results mentioned in 
\ref{Objective.}, our algorithm
does not provide a test for confluence.
We can however use our procedure together with other procedures which
handle \shl-automatic groups,
to prove confluence by an indirect
route, provided the group is \shl-automatic.
Details of the theory of how this is done can be found in
\cite{WPiG}. The practical details are carried out in
programs by Derek Holt---see \cite{Holt:KBMAG}.

We will introduce the concept of \RR-reduction, that is, reduction
using a two-variable automaton, which we call \Rules, encoding our
possibly infinite set of rules.
We prove some results about how reducibility may change with time.

\subsection{Properties of the rule automaton.}
\label{properties of Rules}
The most important data structure is a small two-variable PDFA 
which we call \Rules. Roughly speaking, this accepts all the rules
found so far. It has the following properties.
\begin{enumerate}
\item $\Rules$ is a trim rule automaton.
\item $\Rules$ has one initial state and one final state and they are
equal.
\item $\Rules$ and its reversal
$\Rev \Rules$ are both partially
deterministic.
\item Any arrow labelled $(x,x)$, with either source or target
the initial state, has source equal to target.
has source the initial state. If this condition is not fulfilled, we can
identify the source and target of the appropriate $(x,x)$-arrows,
and then weld. We will still have a rule automaton.
Later on (see Lemmas~\ref{arrows removed 1} and \ref{arrows removed 2}) we
will show that (after any necessary identifications and welding)
we can omit such arrows without loss, and, in fact, with a gain
given by improved computational efficiency.
Apart from the passages proving these lemmas,
we will assume from now on that there are no arrows
labelled $(x,x)$ with source or target the initial state of $\Rules$.
\end{enumerate}

The first three conditions imply that $\Rules$ is welded.
Since $\Rules$ is a rule automaton,
Proposition~\ref{validrules} shows that
each accepted pair $(u,v)\in L(\Rules)$ gives a valid identity $\bar u = \bar v$
in $G$.

\subsection{The automaton $\SLtwo$.}\label{SLtwo definition}
\begin{figure}[htbp]
\begin{center}
{\fontsize{10}{12pt}\selectfont
\begin{picture}(100,140)(0,-20)
\put(0,50){\circle{5}}
\put(-2,39){$1$}
\put(-12.5,50){\vector(1,0){10}}
\put(1.7,51.7){\vector(1,1){30}}
\put(-6,74){$(x,y),$}
\put(-14,64){$x>y$}
\put(1.7,48.3){\vector(1,-1){30}}
\put(-10,27){$(x,y),$}
\put(-2,17){$x<y$}
\put(33.5,83.5){\circle*{5}}
\put(31.5,73){$2$}
\put(33.5,96){\circle{20}}
\put(29.8,103){$<$}
\put(23,112){$(x,y)$}
\put(33.5,16.5){\circle{5}}
\put(31.5,21){$3$}
\put(33.5,4){\circle{20}}
\put(29.8,-8){$<$}
\put(23,-17){$(x,y)$}
\put(67,50){\circle*{5}}
\put(65,39){$4$}
\put(2.5,50){\vector(1,0){62}}
\put(22,54){$(x,\$)$}
\put(35.2,81.8){\vector(1,-1){30}}
\put(52,70){$(x,\$)$}
\put(35.2,18.2){\vector(1,1){30}}
\put(53,25){$(x,\$)$}
\put(67,50){\vector(1,0){62}}
\put(89,54){$(x,\$)$}
\put(131.5,50){\circle*{5}}
\put(129.5,39){$5$}
\end{picture}
}
\end{center}
\caption{\sf The automaton $\SLtwo$. Solid dots represent final states.
Roman letters represent arbitrary
letters from the alphabet $A$ and the labels on the arrows indicate multiple
arrows.
For example, from state $2$ to itself there is one arrow for each pair in
$A\times A$.}

\label{SL2}
\end{figure}
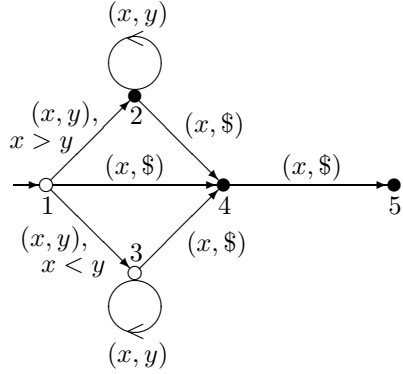

The automaton $\Rules$ may accept pairs $(u,v)$ such that $u$ is shorter
than $v$. We cannot consider such a pair as a rule and so we want to
exclude it. To this end we introduce the automaton $\SLtwo$. This is a five
state automaton, depicted in Figure~\ref{SL2}, which accepts pairs $(u,v)
\in \Astar\times\Astar$, such that $u$ and $v$ have no common prefix, $u$ is 
\shl-greater than $v$ and $|v|
\le |u| \le |v|+2$. By combining $\SLtwo$ with $\Rules$, we obtain a regular
set of rules
$\SetOfRules \Rules$, which is possibly infinite, namely
$L(\Rules)\cap L(\SLtwo)$. An automaton accepting this set can be constructed
as follows. Its states are pairs $(s,t)$, where $s$ is a state of $\Rules$ and
$t$ is a state of $\SLtwo$. Its unique initial state is the pair of initial
states in $\Rules$ and $\SLtwo$. A final state is any state $(s,t)$ such that
both $s$ and $t$ are final states. Its arrows are labelled by $(x,y)$, where
$x\in A$ and $y\in A^+$. Such an arrow corresponds to a pair of arrows, each
labelled with $(x,y)$, the first from $\Rules$ and the second from $\SLtwo$.

\subsection{Restrictions on relative lengths.}
\label{right-hand sides at most two shorter}
The following discussion is closely connected with
\thref{conditions and minimals}.
The restriction $|u| \le |v| +2$ needs some explanation. The point is that if
we have a rule with $|u| > |v| +2$, then we have an equality $\bar u = \bar v$
in $G$. We write $u = u'x$, where $x\in A$. The formal inverse $X$ of $x$
is also an element of $A$. We therefore have a pair of words
$(u', vX)$ which represent equal elements in $G$. If our set of rules were
to contain such a rule, then $u=u'x$ would reduce to $vXx$, and this reduces
to $v$, making the rule $(u,v)$ redundant. This leads to an obvious technique
for transforming any rule we find into a new and better
rule with $|v| \le |u| \le |v|+2$. Since we take this into account when
constructing the automaton $\Rules$, we are justified in making the
restriction.

This analysis can be carried further. Let
$u=u_1\cdots u_{r+2}=u'u_{r+2} = u_1 u''$ and
let $v=v_1\cdots v_r$. If $u_1>v_1$, then the rule $(u,v)$ can be
replaced by the better rule $(u',vu_{r+2}^{-1})$.
If $u_2 > u_1^{-1}$, then $(u,v)$ can be replaced by $(u'',u_1^{-1}v)$.
We do in fact carry out these steps when installing new rules.
The extra information could have been included in the FSA $\SLtwo$. However,
it seems that this would involve more complicated coding at various 
points, probably without any gain in efficiency.

We could consider the steps just 
described as an attempt to force our structures to define a set of 
rules which conforms to known properties (see \thref{conditions and
minimals}) of the set of
$U$-minimal rules (see \ref{U} for the definition of $U$).
The most important reason for insisting on these additional 
restrictions on our rules is to keep down the size of our data 
structures.

\subsection{The basic structures.}
The basic structures used in our procedure are:
\begin{enumerate}
\item A two-variable automaton $\Rules$ satisfying the conditions
laid down in \ref{properties of Rules}.
When we want to
specify that we are working with the $\Rules$ automaton during the $n$th 
Knuth--Bendix pass (see \ref{Knuth--Bendix pass} for the definition of a
Knuth--Bendix pass), we will use the notation $\Rules[n]$. We extract
explicit rules from $\Rules[n]$ by taking elements of the intersection
$\SetOfRules{\Rules[n]} = L(\Rules[n]) \cap L(\SLtwo)$. The two-variable
automaton $\SLtwo$ was defined in Section~\ref {SLtwo definition} and
is depicted in Figure~\ref{SL2}. 

\item A finite set \Store of rules, which is the disjoint union of
several subsets of rules : \Considered, \This, \New and \Delete.
One point of the separate subsets is to avoid constantly doing the 
same critical pair analyses. Another point is to ensure that our 
Knuth--Bendix process is fair (see \thref{fair}). The reason for 
holding some rules in a \Delete list, rather than delete them
immediately, is to make reduction more efficient.
This will be explained further in \ref{deletion delay}.

\Store will continually change, while $Rules$ is constant during a
Knuth--Bendix pass. We change $\Rules$
at the end of each Knuth--Bendix pass.
We will perform the Knuth--Bendix process, using the rules in \Store
for critical pair analysis, as described in \ref{Critical pair analysis.}.
\item \Considered is a subset of \Store such that each rule
has already been compared with each other rule in \Considered,
including with itself,
to see whether left-hand sides overlap.
The consequent critical pair
analysis has also been carried out for pairs of rules in \Considered.
Such rules do not need to be compared with
each other again.
\item \This is a subset of \Store (empty at the beginning of each
Knuth--Bendix pass) containing rules which we plan to use
during this pass to compare for overlaps with the rules in \Considered, as in
\ref{first case}. These rules are minimal for
the current pass (see \ref{minimization}) and so should not be minimized again.
\item \New is a subset of \Store containing new rules which have
been found during the current pass, other than those which are output
by the minimization routine (see \ref{minimization} for the meaning of
``minimization''). Rules which are
output by the minimization routine are added to \This.
\item \Delete is a subset of \Store containing rules which are to be
deleted at the end of this pass.
\item The two-variable automaton \WDiff contains
all the states and arrows of $\Rules[n]$, and possibly other states and
arrows. It satisfies the conditions of \ref{properties of Rules}.
This automaton is used to accumulate appropriate new rules which are 
output by the minimization routine. As rules are considered during the 
Knuth--Bendix pass, states and arrows of \WDiff are marked as \needed.
At the end of the pass, other states and arrows are removed, and 
\WDiff becomes the new \Rules automaton $\Rules[n+1]$.
\item A PDFA $P(\Rules)$ formed from
\Rules by a certain subset construction.
This automaton accepts words which are \RR-reducible, that is, 
words which contain a left-hand side of a rule in
 $\SetOfRules\Rules$. The automaton is used as part of 
our rapid reduction procedure (see \thref{Fast reduction}).
More details of $P(\Rules)$ are provided in \ref{P}.
\item A PDFA $Q(\Rules)$ which accepts the reversals of left-hand sides 
of rules in  $\SetOfRules\Rules$. This is also formed from \Rules by 
a subset construction and is also used for rapid reduction.
 More details of $Q(\Rules)$ are 
provided in \ref{Q}.
\end{enumerate}

\subsection{Initial arrangements.}
Before describing the main Knuth--Bendix process, we explain how the
data structures are initially set up. Let $R$ be the original set of
defining relations together with special
rules of the form $(x.\inverse(x),\epsilon)$ which
make the formal inverse $\inverse(x)$ into the actual inverse of $x$.

We rewrite each relation of $R$ in the form
of a relator, which we cyclically reduce in the free group.
We assume that each relator has the form
$l.\inverse(r)$, where $l$ and $r$ are elements of $\Astar$ and $(l,r)$ is
accepted by $\SLtwo$.

For each rule $(l,r)$, including the special rules
$(x.\inverse(x),\epsilon)$, we form a rule
automaton, as explained in \ref{ruleworddiff}. These automata are then welded
together to form the two-variable rule automaton \WDiff satisfying the 
conditions of \ref{properties of Rules}. Each state and arrow 
of \WDiff is marked as \needed.
Each of these rules is inserted into \New. \Considered, \This and 
\Delete are initially empty. Set $\Rules[1]=\WDiff$.

\subsection{The main loop---a Knuth--Bendix pass.}
\label{main loop}
We now describe the procedure followed during the course of a single
Knuth--Bendix pass.

A significant proportion of the time in a Knuth--Bendix pass is spent in
applying a procedure which we term {\it minimization}. Each rule encountered
during the pass is input (often after a delay)
to this procedure and the output is called a 
{\it minimal rule}. 
The details of this process are given in sections~\ref{minimization}
and \ref{handling minimization output}.

\begin{enumerate}
\item\label{start main loop}
At the beginning of a Knuth--Bendix pass, \This is empty.
If $n>0$, save space by deleting previously defined automata $P(\Rules[n])$,
$Q(\Rules[n])$ and $\Rules[n]$. Increment $n$.
The integer $n$ records which Knuth--Bendix pass we are currently working on.
\item\label{process Considered}[Step]
For each rule $(\lambda,\rho)$ in \Considered,
minimize $(\lambda,\rho)$ as in \ref{minimization} and handle the output rule
$(\lambda_1,\rho_1)$ as in \ref{handling minimization output}.
This may affect \Store and \WDiff.
\item \label{process New}[Step]
For each rule $(\lambda,\rho)$ in \New,
minimize $(\lambda,\rho)$ as in \ref{minimization} and handle the output as in
	\ref{handling minimization output}.
This may affect \Store and \WDiff.

Since rules added to \New during minimization are always strictly smaller than
the rule being minimized (see \ref{inserting}), it follows that
the process of examining rules in \New does not continue indefinitely.
As a result, we can be sure that our process is fair (see \ref{fair}).
\item\label{process This}
For each rule $(\lambda,\rho)$ in \This:
	\begin{enumerate}
	\item Delete the rule from \This and add it to \Considered.
	\item\label{compare for overlaps}[Step]
		 For each rule $(\lambda_1,\rho_1)$ in \Considered:
		\begin{quote} Look for overlaps between $\lambda$ and
		$\lambda_1$. That is we have to find each
		suffix of $\lambda$ which is a prefix
		of $\lambda_1$ and
		each suffix of $\lambda_1$ which is a
		prefix of $\lambda$.
		Then \RR-reduce in two different ways as in
		\ref{first case}, obtaining a pair of words $(u,v)$
		with $u\ge v$. (Roughly speaking, \RR-reduction means
		the use of rules in $\SetOfRules{\Rules}$. More
		precision is provided in \ref{inserting}.)
		If $u>v$, $(u,v)$ is inserted into \New,
		 unless it is already in \Store.

		Note that we may have to allow $\lambda=\lambda_1$ in order to
		deal with the case where two different rules have the
		same left-hand side. In this case, both the prefix and suffix 
		of both left-hand sides is equal to $\lambda=\lambda_1$.
		\end{quote}
	\end{enumerate}
\item\label{using WDiff}
\WDiff was possibly affected in \thref{process Considered} and
\thref{process New}.
With \WDiff in its present form,
delete from \WDiff all arrows and states which are not marked as
\needed. Copy \WDiff into $\Rules[n+1]$
and mark all arrows and states of \WDiff as not \needed.
\item \label{process Delete}
Delete the rules in \Delete.
\item This ends the description of a Knuth--Bendix pass.
Now we decide whether to terminate the Knuth--Bendix process. Since we know of 
no procedure to decide confluence of an infinite system of rules (indeed, it is
probably undecidable), this decision is taken on heuristic grounds. In our
context, a decision to terminate could be taken simply on the grounds that 
\WDiff and $\Rules[n]$ have the same states and arrows. In other words, no new
word-differences or arrows between word-differences have been found or
deleted during this pass.
If the Knuth--Bendix process is not terminated, go to \ref{start main loop}.
\end{enumerate}

\begin{definition}
\label{minimization}[Definition]
We now provide the details of the minimization routine. 
This processes a rule so as to create from it a minimal rule
(see \thref{minimal}), where, roughly speaking, minimality is defined using
the current set of rules.
Since the set of rules is changing, this is a bit difficult to pin
down. So instead we make the following definition, which is more
precise, though the underlying concept is the same.
Let $(u,v)\in \Astar\times\Astar$ and let $u=u_1\cdots u_p$
and $v=v_1\cdots v_q$, where $u_i, v_j \in A$.
We say that $(u,v)$ is a \textit{minimal rule} if $u\neq v$,
$\bar u = \bar v$ in $G$
and the following procedure does not change $(u,v)$. The procedure is
called the \textit{minimization routine}.
We always start the minimization routine with $u>v$, though this
condition is not necessarily maintained as $u$ and $v$ change during the
routine.
Here the meaning of a ``minimal rule'' changes with time: a rule may
be minimal at one time and no longer minimal at a later time.

\begin{enumerate}
\item\label{reduce prefix}
\RR-reduce (that is, reduce using the rules of \Rules)
the maximal proper prefix $u_1\cdots u_{p-1}$
of $u$ obtaining $u'$.
Reduction may result in rules being added to \New as described in
\ref{reduction gives new rule}.
If $u\neq u'u_p$, change $u$ to $u'u_p$ and go to
Step~\ref{reduce u}.
\item\label{reduce suffix}
\RR-reduce the maximal proper suffix $u_2\cdots u_p$
of $u$ obtaining $u''$.
Reduction may result in new rules being added to \New.
Replace $u$ by $u_1u''$.
\item\label{reduce u}
If $u$ has changed since the original input to the minimization routine,
then \RR-reduce $u$ as explained in \ref{history stack}.
This may result in rules being added to \New as described in
\ref{reduction gives new rule}.
\item\label{reduce v}[Step]
\label{minimization loop}[Step]
\RR-reduce $v$.
\item If $v > u$, interchange $u$ and $v$.
\item If (a) $p >q+2$ or (b) if $p = q+2$, $q>0$
and $u_1 > v_1$ or (c) if $p=2$, $q=0$ and
$u_1>\inverse(u_2)$,  replace $(u,v)$ by
$(u_1\cdots u_{p-1}, v_1\cdots v_q \inverse(u_p))$ and repeat this step until
we can go no further.
\item If $p=q+2$ and $u_2 > \inverse(u_1)$, replace $(u,v)$ by
$(u_2\cdots u_p,\inverse(u_1)v_1\cdots v_q)$.
\item If $q>0$ and $u_1 = v_1$, cancel the first letter
from $u$ and from $v$ and repeat this step.
\item If $q>0$ and $u_p =v_q$, cancel the last letter from $u$ and
from $v$ and repeat this step.
\item If $(u,v)$ has changed since the last time Step~\ref{minimization
loop} was executed, go to Step~\ref{minimization loop}.
\item Output $(u,v)$ and stop.
\end{enumerate}
\end{definition}

Note that the output could be $(\epsilon,\epsilon)$, which means
that the rule is redundant. Otherwise we have output $(u,v)$ with $u>v$.
Note that the minimization procedure keeps on decreasing $(u,v)$ in the
ordering given by using first the \shl-ordering on $u$ and then, in
case of a tie, the \shl-ordering on $v$.
Since this is a well-ordering, the minimization
procedure has to stop.

\subsection{Handling minimization output.}
\label{handling minimization output}
Suppose the input to minimization is $(\lambda,\rho)$ and its output
is $(\lambda_1,\rho_1)$.
\begin{enumerate}
\item If $(\lambda_1,\rho_1) \neq (\epsilon,\epsilon)$,
incorporate (by welding) $(\lambda_1,\rho_1)$ into the language accepted
by \WDiff.
Insert $(\lambda_1,\rho_1)$ into \This if it was not already in \This
or \Considered. Remove it from \New, if it was there previously.
\item\label{deletion: proper subword}
If some proper subword of $\lambda$ is \RR-reducible, then this will
be discovered during the first few steps of minimization.
($(\lambda_1,\rho_1) = (\epsilon,\epsilon)$ turns out to be
a special case of this,
as we will see in \ref{deleting rules irred}.)
In this case, delete $(\lambda,\rho)$ from \Store immediately the
minimization procedure is otherwise complete.
\item\label{deletion delay}
If, at the time of minimization,
all proper subwords of $\lambda$ were \RR-irreducible and
if $(\lambda,\rho)$ was not minimal, move
$(\lambda,\rho)$ to the \Delete list. The reason for this possibly
surprising policy of not deleting immediately is that further
reduction during this pass may once again
produce $\lambda$ as a left-hand side by the
methods of \ref{Fast reduction} and
\ref{finding the left-hand side}. We want to avoid the work involved in
finding the right-hand side by the method which will be explained in
\ref{finding the right-hand side}. For this, we need to have a rule
in \Store with left-hand side equal to $\lambda$---see \ref{reduction
gives new rule}.
\end{enumerate}

\subsection{Details on the structure of \WDiff.}
\label{WDiffdetails}
At the beginning of Step~\ref{using WDiff}, each state $s$ of $\WDiff$
is associated to a word $w_s \in A\uast$ which is irreducible with respect
to $\SetOfRules{\Rules[n]}$. $\WDiff$ is a rule automaton: the rule
automaton structure is given by associating the
element $\overline{w_s}\in G$ to the state $s$.
Whenever a minimal rule $r$ is encountered during the $n$th pass, it 
is adjoined to the accepted language of \WDiff by welding and the
corresponding states and arrows are marked as \needed.
State labels are calculated as and when new states and arrows are added 
to $\WDiff$.

At the end of the $n$th Knuth--Bendix pass,
\WDiff is an automaton which represents the
word-differences and arrows between them encountered during that pass.
At this stage the word
attached to each state is irreducible with respect to the rules in
$\SetOfRules{\Rules[n]}$ but not necessarily 
with respect to the rules implicitly contained in \WDiff. Before starting the
next pass, we \RR-reduce the state labels of \WDiff with respect to 
$\SetOfRules{\WDiff}$. If \WDiff now contains distinct states labelled by the
same word we connect them by epsilon arrows and replace \WDiff by
$\Weld \WDiff$. We then repeat this procedure until
all states are labelled
by distinct words which are irreducible with respect to 
$\SetOfRules{\WDiff}$. If during this procedure a state or arrow marked
as \needed is identified with another which may or may not be marked
as \needed, the resulting state or arrow is marked as \needed.

\subsection{\RR-reduction and inserting rules.}\label{inserting}
Given a word $w$, we look for an \RR-reducible subword $\lambda$ such
that all proper subwords of $\lambda$ are \RR-irreducible,
by looking in $\SetOfRules{\Rules}$.
Later (\thref{Fast reduction})
we will describe how to
do this quickly, but, at the moment, the reader can just think of a
non-deterministic search in the automaton giving the \shl rules
recognized by $\Rules$. Having found a reducible subword
$\lambda$ of $w$, with no reducible subword,
we do not automatically use the corresponding right-hand
side $\rho$, found from the exploration of \Rules,
because this naive approach is computationally inefficient.
Instead we look in \Store to see if there is a rule
$(\lambda,\rho)$. If there is such a rule, then we can find it
quickly given $\lambda$,
and we proceed with our reduction, replacing the subword
$\lambda$ in $w$ with $\rho$.

It may however turn out that
we can find an \RR-reducible subword $\lambda$ of $w$, with no
\RR-reducible subwords, and yet
there is no rule of the form $(\lambda,\rho)$ in \Store.
In this case, we have to spend time
finding such a rule in $\SetOfRules\Rules$.
Once found, we immediately insert it into \Store,
otherwise the logic of the Knuth--Bendix procedure can go wrong.

In this way, reduction of a single word can result in the insertion of
several new rules into \Store.

It follows from the above description
that the \RR-reducibility of a word $w$ depends only on \Rules.
Since \Rules does not change during a Knuth--Bendix pass,
exactly the same subset of $\Astar$ will be \RR-reducible throughout 
such a pass. However, because we may use rules in the changing set \Store,
the \textit{result} of \RR-reduction may change during a pass.

Another, more conventional, source of rules to insert into \Store come
from critical pair analysis in \thref{compare for overlaps}.

Minimization also results in rules being added to \Store, both
directly, as the output of the minimization procedure, but also
indirectly because minimization uses reduction, and, as we will see in
\ref{finding the right-hand side}.  reduction can add rules to \Store.
It is important to note that any rules
added to \Store during the minimization of a rule $(\lambda,\rho)$ are
strictly smaller than $(\lambda,\rho)$, if we order such pairs by
using $\lambda$ first and then $\rho$ in case of a tie.
We used this fact when discussing \thref{process New}.

\subsection{Deleting rules.}\label{deleting rules}
Deletion of rules happens only at the end of each minimization step, and at the
end of each pass, when rules marked for deletion are actually deleted.
During a Knuth--Bendix pass, deletion
does not occur after the beginning of Step~\ref{process This}.
Suppose that the output from minimization of $(\lambda,\rho) \in
\Store$ is $(\lambda_1,\rho_1)$.

\begin{enumerate}
\item\label{deleting rules irred}[Case]
If every proper subword of $\lambda$ is \RR-irreducible, then
$\lambda_1$ is a non-trivial subword of $\lambda$.
This follows by going through the successive steps of
minimization (\thref{minimization}).
These change $\lambda$ and $\rho$, while maintaining
the inequality $\lambda > \rho$.
In particular $\lambda_1 > \rho_1$, so that $\lambda_1 \neq \epsilon$.
If $(\lambda_1,\rho_1) \neq (\lambda,\rho)$, then we delete
$(\lambda,\rho)$ after a delay.
The mechanism is to mark it for
deletion by moving it to the \Delete list and actually delete it only
at the end of the current Knuth--Bendix pass (Step~\ref{process Delete}).

\item\label{deleting rules red}[Case]
If some proper subword of $\lambda$ is reducible,
then $(\lambda,\rho)$ is immediately deleted from \Store
at Step \ref{deletion: proper subword} at the end of the
minimization procedure. (\RR-reducibility of some proper subword of
$\lambda$ is discovered
at Step~\ref{reduce prefix} or \ref{reduce suffix}.)
\end{enumerate}

\begin{lemma}\label{before This}
Suppose that, for some $n\in\naturals$,
there is a rule $(\alpha,\beta)\in\Store$ 
during the $n$-th Knuth--Bendix pass,
{\em before} the beginning of Step~\ref{process This}.
Then there is a non-trivial subword $\lambda$ of
$\alpha$ such that some rule
$(\lambda,\rho)$ is output from some instance of the minimization
procedure during the $n$-th pass.
If $\lambda=\alpha$, then $\rho \le \beta$.
The rule $(\lambda,\rho)$ is a rule in \Store at the
beginning of the $(n+1)$-st pass and is accepted by $\Rules[n+1]$.
\end{lemma}
\begin{proof}
By examining \ref{main loop}, we see that
$(\alpha,\beta)$ must be the input to the minimization routine
at some time during the $n$-th
pass. (We check the four possibilities, namely that it is in
\Considered, \This, \New or \Delete, one by one. If it is in \Delete,
it must have been the input to the minimization procedure at some
earlier stage during the $n$-th pass.)

We first deal with the case where
some proper subword of $\alpha$ is \RR-reducible during the
$n$-th pass.
During the first three steps of minimization (\thref{minimization}),
an \RR-reducible subword $\lambda$ of $\alpha$ is found, with the
property that all the proper subwords of $\lambda$ are \RR-irreducible.
Minimization then either finds a rule of the form
$(\lambda,\rho)$ already in \Store, or such a rule is
added to \New by the reduction process---see \ref{reduction gives new
rule}.
In any case, it will either be minimized during this pass, or it has
already been minimized (and possibly moved to the \Delete list.

At the moment when $(\lambda,\rho)$ is minimized during the $n$-th pass,
we must be in Case~\ref{deleting rules irred}.
So the output $(\lambda_1,\rho_1)$ from the minimization
procedure with input $(\lambda,\rho)$ gives the required rule.
$\lambda_1$ is a subword of $\lambda$ and $\lambda$ is a proper
subword of $\alpha$.

Alternatively, all proper subwords of $\alpha$ are \RR-irreducible
during the $n$-th pass,
in which case we set $(\lambda,\rho)$ to be the output from 
minimization of $(\alpha,\beta)$.
By \ref{deleting rules irred}, $\lambda$ is a non-trivial subword 
of $\alpha$. If $\lambda=\alpha$, then $\rho\le \beta$.
\end{proof}

\begin{lemma}\label{after This}
Suppose that, for some $n\in\naturals$,
there is a rule $(\alpha,\beta)\in\Store$ 
during the $n$-th Knuth--Bendix pass, {\em after} the beginning of
Step~\ref{process This}.
Then there is a non-trivial subword $\lambda$ of $\alpha$ such that some rule
$(\lambda,\rho)$ is output from some instance of the minimization
procedure during the $(n+1)$-st pass.
If $\lambda=\alpha$, then $\rho \le \beta$.
\end{lemma}

\begin{proof}
If $(\alpha,\beta)$ is in the \Delete list, then it must have been
input to the minimization procedure at some earlier time during the
$n$-th pass. By \thref{deleting rules red}, every proper subword of
$\alpha$ must have been found to be \RR-irreducible during the $n$-th pass.
Let $(\alpha',\beta')$ be the output from minimization.
By \thref{deleting rules irred}, $\alpha'$ is a non-trivial subword
of $\alpha$, and, if $\alpha'=\alpha$, then $\beta'<\beta$.
Now $(\alpha',\beta')$ is in \Store at the beginning of the $(n+1)$-st
pass. We apply \thref{before This} to $(\alpha',\beta')$ at the
$(n+1)$-st pass.

If $(\alpha,\beta)$ is not on the \Delete list, then it must be in \Store
at the beginning of the $(n+1)$-st pass.
Once again, we can apply \thref{before This}.
\end{proof}

The following result is often applied with $w=\alpha$.
\begin{proposition}\label{stays RR reducible}
Let $w\in \Astar$ be a word which contains the left-hand side $\alpha$ of a
rule $(\alpha,\beta)$ input to the minimization routine during the
$n$-th Knuth--Bendix pass.
Then, for $m\ge n$, $w$ contains the left-hand side of a rule which is
input to the minimization procedure during the $m$-th Knuth--Bendix pass.
Moreover $w$ is \RR-reducible for $m>n$.
\end{proposition}
\begin{proof}
We assume inductively that if $m>n$ then
$w$ contains a subword $\alpha$, such that a rule of the form $(\alpha,\beta)$
is input to the minimization procedure during the $(m-1)$-st pass.
Since minimization happens only
before the beginning of Step~\ref{process
This}, \thref{before This} gives a rule $(\lambda,\rho)$, such that
$\lambda$ is a non-trivial subword of $\alpha$.
Moreover, $(\lambda,\rho)$ is
minimal during the $(m-1)$-st pass and is contained in \Store at the
beginning of the $m$-th pass. Therefore $(\lambda,\rho)$ is input to the
minimization procedure during the $m$-th pass, as required.

The rule $(\lambda,\rho)$ is welded into \WDiff during the $(m-1)$-st
pass and is therefore accepted by $\Rules[m]$.
It follows that $w$ is \RR-reducible during the $m$-th pass.
Inductively this is true for all $m>n$.
\end{proof}

%% file: Sec6.tex
\section{Correctness of our Knuth--Bendix Procedure}
\label{correctness}
In this section we will prove that the procedure set out in
Section~\ref{our version} does what we expect it to do.
One hazard in programming Knuth--Bendix is that some seemingly
clever manoeuvre changes the Thue equivalence relation.
The key result here is \thref{maintain congruence}, which carefully
analyzes the effect of our various operations on Thue equivalence. In fact
it provides more precise control, enabling other hazards, such as
continual deletion and re-insertion of the same rule, to be avoided.
It is also the most important step in proving our main result,
\thref{main}. This says that if our program is applied to a group
defined by a regular set of minimal rules, then, given sufficient time
and space, a finite state automaton accepting exactly these rules will
eventually be constructed by our program, after which the program will loop
indefinitely, repeatedly reproducing the same finite state automaton (but
requiring a steadily increasing amount of space for redundant information).

\begin{definition}\label{time}[Definition]
For a discrete time $t$, we denote by \Storen t the rules in
\Store at time $t$ in our Knuth---Bendix procedure. We take $t$ to be
the number of elementary steps since the start of the program,
assuming the program is expressed in some sort of pseudocode.
Any other similar measure of time would do equally well.
\end{definition}

\begin{definition}
A quintuple $(t,s_1,s_2,\lambda,\rho)$, where $t$ is a time, and $s_1$,
$s_2$, $\lambda$ and $\rho$ are elements of $\Astar$, is called an
{\it elementary \Storen t-reduction} $u \rightarrow_{\Storen t} v$ from $u$
to $v$ if $(\lambda,\rho)$ is a rule in \Storen t, $u=s_1\lambda s_2$
and $v=s_1\rho s_2$. We call $(\lambda,\rho)$ the \textit{rule associated
to the elementary reduction}.
\end{definition}

We now define the main technical tool that we will use in this section.
\begin{definition}\label{Thue path}
Let $t\ge 0$.
By a \textit{time-$t$ Thue path} between two words $w_1$ and $w_2$,
we mean a finite
sequence of elementary \Storen t-reductions and inverses of elementary
\Storen t-reductions connecting $w_1$ to $w_2$, such that none of the
rules associated to the elementary reductions is in \Delete at time $t$.
We talk of the words which are the source or target of these elementary
reductions as \textit{nodes}.
The path is considered as having a direction from $w_1$ to $w_2$.
The elementary reductions in our path will be consistent with this
direction and will be called \textit{rightward} elementary reductions.
The inverses of elementary reductions in our path
will be in the opposite direction
and will be called \textit{leftward} elementary reductions.
\end{definition}

All our insertions and deletions of rules have been organized so that the
following result holds.
\begin{proposition}
\label{consistent rules}
Let $\langle A / R \rangle$ be the finite presentation of a group $G$
at the start of the Knuth--Bendix process.
Then the group defined by subjecting the free group generated
by $A$ to all relations of the form $\lambda =\rho$ as
$(\lambda,\rho)$ varies over \Storen t is at all times $t$ isomorphic to $G$
with the isomorphism being induced by the unchanging map $A\to G$.
\end{proposition}

\begin{proposition}
\label{maintain congruence}
Let $t\geq0$ and suppose that we have a Thue path from $u$ to 
$v$ in \Storen t with maximum node $w$. Then for any time $s\geq t$,
there exists a time-$s$ Thue path from $u$ to $v$
with each node less than or equal to $w$.
\end{proposition}
\begin{proof}
Note that, given a Thue path, we may assume, if we wish,
that no node is repeated, because we could shorten the path to avoid
repetition.
We show by induction on $s$ that, if at some time $t\leq s$ there is a Thue
path between words $u$ and $v$ with all nodes no bigger than $\max(u,v)$, then
there is also such a Thue path at time $s$. So suppose that we have proved this
statement for all times $s'<s$. 

We first consider the special case where $r_0=(u,v)$ is
a rule being input to the minimization routine (see Definition~\ref{minimization}) 
at time $t$, and $s$ is the time at the
end of the subsequent invocation of the
minimization handling routine \ref{handling minimization output}. 
There is a Thue path (of length one) from $u$ to $v$ at time $t$.
By induction we are assuming that at time $s-1$
there is a Thue path from $u$ to $v$ with maximum node $u$.
We must show that there is such a Thue path at time $s$.

One possibility is that
$r_0$ is already minimal, in which case there is a Thue path of length
one from $u$ to $v$, both at the beginning and at the end of
minimization. So we assume that $r_0$ is not minimal.
Then the last step in \ref{handling minimization output}
is that either $r_0$ is placed in the \Delete list or else $r_0$ is
simply deleted immediately.

What we need to show therefore is that the Thue path $p$ from $u$ to
$v$, which exists at time $s-1$,
does not use an elementary reduction coming from $r_0$. It is
part of our inductive hypothesis that the largest node occurring on $p$ is
$u$, and we have already pointed out that we can assume there is no
repetition of nodes along $p$.

Each step of minimization takes an input pair of words and outputs
a possibly different pair of words which is used as the input to the
next step.
The initial input is $r_0=(u,v)$ and the final output is
either $r_n=(\epsilon,\epsilon)$ or a minimal rule $r_n= (u',v')$.
Let $r_0,r_1,r_2,\ldots,r_n$ be the sequence of such inputs and outputs in the
minimization of $(u,v)$.
By considering each step of minimization in turn,
we will show that for each $i$, $1\leq i\leq n$, if there is a time-$s$
Thue path between the two sides of $r_i$ with maximum node no bigger
than either side of $r_i$, then there is a time-$s$ Thue path between
the two sides of $r_{i-1}$ with maximum node no bigger than either side
of $r_{i-1}$.
We then obtain the desired time-$s$ Thue path between
$u$ and $v$ by using descending induction on $i$.
This is a subsidiary induction to our main induction on $s$.
The base case $i=n$
is true, since at time $s$ the rule $r_n$ has been installed in \Store.

To make the task of checking the proof easier, we use the same numbering
and notation here as in Definition~\ref{minimization}.
\begin{enumerate}
\item At the end of the current step,
there is a sequence of elementary reductions
from $u_1\ldots u_{p-1}$ to $u'$, but this may not constitute a Thue path
since some of the associated rules may be in \Delete. However, any such rule 
$(\lambda,\rho)$ in \Delete will, at some time $s' < s$,
have been in \Store but not in \Delete.
Therefore, by our induction on $s$, at time $s-1$ there is
a Thue path $p$ from $\lambda$ to $\rho$ with maximum node $\lambda$.
Now $\lambda \le u_1\ldots u_{p-1} < u$ and so $\lambda$ is smaller than the
left-hand side of $r_0$. Therefore $r_0$ cannot be used in $p$.
So $p$ continues to be a Thue path at time $s$.
This completes the downward induction step on $i$ in this case.

\item This step is analogous to the previous step.

\item The sequence of \RR-reductions of $u$ to the 
current left-hand side does not use the rule $r_0$
and so the required Thue path exists by induction on $s$. 

\item Let $v'$ be the \RR-reduction of $v$. Immediately after this step there
is a Thue path from $v$ to $v'$ with maximum node $v$ which does not use $r_0$.
By the induction hypothesis on $s$, there is such a Thue path at time
$s-1$.
Since it does not use $r_0$, it continues to be a Thue path at time $s$.
Hence a time-$s$
Thue path from $u$ to $v'$ with maximum node either $u$ or $v'$ yields a 
time-$s$ Thue path from $u$ to $v$ with maximum node $u$ or $v$.
(Recall that, because of previous steps which may shorten $u$,
$u$ may be smaller than $v$ at this point.)
This completes the downward induction step on $i$ in this case.

\item If there is a Thue path from $u$ to $v$ with maximum node either
$u$ or $v$, then the reverse of this path is a Thue path from $v$ to $u$.

\item Suppose that the input to this step is $(u'x,v)$.
Then the output is either the same as the input or is equal to
$(u',v.\inverse(x))$, with $u' > v.\inverse(x)$.
In the first case there is nothing to prove. In the latter 
case, we have by our downward induction on $i$
a time-$s$ Thue path from $u'$ to 
$v.\inverse(x)$ with maximum node $u'$.
This will give a time-$s$ Thue path from $u'x$
to $v.\inverse(x)x$ with maximum node $u'x$.
Furthermore, at the beginning of the Knuth--Bendix process, there was
a Thue path of length one
from $\inverse(x)x$ to $\epsilon$ with maximum node equal
to $\inverse(x)x$. Therefore, by our
induction hypothesis, there is such a path at time $s-1$, just before
possible deletion of $r_0$. Now $u'x > v.\inverse(x)x \ge
\inverse(x)x$. So the time-$(s-1)$ Thue path from $\inverse(x)x$ to $\epsilon$
cannot use $r_0$, and it
remains a Thue path at time $s$. It follows that
there is a Thue path from $u'x$ to $v$ with maximum node $u'x$ at time $s$.

\item This step is analogous to the previous step. 

\item If the input to this step is $(xu',xv')$ then the output is $(u',v')$.
A time-$s$ Thue path from $u'$ to $v'$ with maximum node $u'$ yields
a time-$s$ Thue path from $xu'$ to $xv'$ with maximum node $xu'$.

\item This step is analogous to the previous step. 

\end{enumerate}

This completes the induction on $s$ for the special case
where $r_0=(u,v)$ is
a rule being input to the minimization routine
(see Definition~\ref{minimization}) 
at time $t$, and $s$ is the time at the
end of the subsequent invocation of the
minimization handling routine \ref{handling minimization output}. 
Now consider the general case, again assuming the induction statement
true at time $s-1$.
The only reason why a Thue path at time $s-1$ between $u$ and
$v$ will not work at time $s$ is if some elementary reduction used in this
path has an associated rule $(\lambda,\rho)$ in \Storen {s-1} which is
deleted at time $s$.
Since deletion only takes place as a result of minimization, we know
that what must be happening is that we are right at the end of
minimizing $(\lambda,\rho)$, with minimization completing exactly at time $s$.
But the special case already proved shows
that there is a time-$s$ Thue path between $\lambda$ and $\rho$ with
no node bigger than $\lambda$. Therefore the time-$(s-1)$
Thue path can always be
replaced by a time-$s$ Thue path without increasing the maximum node. 
\end{proof}

\begin{lemma}\label{stays S-reducible}
If a word is \Storen t-reducible, it is \Storen s-reducible for all
$s>t$.
\end{lemma}
\begin{proof}
If $u$ is \Storen t-reducible, there is an elementary
\Storen t-reduction $u\to_{\Storen t} v$.
This means that $v<u$. By Proposition~\ref{maintain congruence},
for each time $s>t$, there is a Thue path from $u$ to $v$ with maximum
node $u$. The first elementary reduction in this path has the form
$u\to w$ at time $s$. This proves the result.
\end{proof}

\begin{lemma}\label{no repeats}
At any time $t$, $\Storen t$ is a list of rules which contains no duplicates.
If a rule is deleted from \Store, it will never be re-inserted.
(Here we mean actual deletion, not just placing the rule on the
\Delete list for future deletion.)
\end{lemma}
\begin{proof}
The first statement follows by looking through \ref{main loop} and
checking where insertions of rules take place. We always take care not
to insert a rule a second time if it is already present.

Let $(\alpha,\beta)$ be a rule which is deleted at time $s$.
We assume by contradiction that it is re-inserted at a later time $t$.
We choose $m$ and $n$ so that time $s$ occurs during the $m$-th
Knuth--Bendix pass and time $t$ during the $n$-th. Then $m \le n$.

We note that all proper subwords of $\alpha$ are \RR-irreducible
during the $m$-th pass. For otherwise \thref{stays RR reducible} shows
that $\alpha$ is \RR-reducible during the $n$-th pass.
But no rule with left-hand side $\alpha$ could then be introduced
during the $n$-th pass, a contradiction.

It follows that we are in Case~\ref{deleting rules irred}.
Therefore $(\alpha,\beta)$ was input to the minimization procedure
during the $m$-th pass and was then moved to \Delete. The actual
deletion took place at the end of the $m$-th pass. It follows that
$n>m$.
The output from the minimization procedure was a rule
$(\lambda,\rho)$, where $\lambda$ is a subword of $\alpha$.
The rule $(\lambda,\rho)$ is welded into \WDiff and is accepted by
$\Rules[m+1]$.
As in the preceding paragraph, we see that $\lambda$ cannot be a
proper subword of $\alpha$, and so $\lambda = \alpha$ and
$\rho<\beta$.
We write $\beta_{m-1} = \beta$ and $\beta_m = \rho$.

Proceeding in this way,
we see that between times $s$ and $t$, rules
of the form $(\alpha,\beta_{i-1})$ ($m \le i \le n$)
are input to the minimization
procedure during the $i$-th Knuth--Bendix pass,
with output $(\alpha,\beta_i)$ where $\beta_i \le \beta_{i-1}$ and
$\beta_m < \beta_{m-1}$.
The rule $(\alpha,\beta_i)$ is produced during the $i$-th
Knuth--Bendix pass and is accepted by $\Rules[i+1]$ for $m\le i \le n$.

It follows that $\alpha$ is \RR-reducible during the $n$-th pass.
Therefore no rule with left-hand side $\alpha$ could be introduced
into \Store as a result of critical pair analysis.
We see from \ref{inserting} that
any rule with left-hand side equal to $\alpha$ which is introduced into
\Store as a result of \RR-reduction during the $n$-th pass
must be of the form $(\alpha,\gamma)$, where $\gamma \le \beta_n<\beta$.
This completes the proof of the contradiction.
\end{proof}

\begin{definition} We say that a word $u$ is \textit{permanently
irreducible} if there are arbitrarily large times $t$ for which
$u$ is \Storen t-irreducible. By Lemma~\ref{stays S-reducible}
this is equivalent to saying
that $u$ is \Storen t-irreducible at all times $t\ge 0$.
A rule $(\lambda,\rho)$ in \Store
is said to be \textit{permanent} if $\rho$ and
every proper subword of $\lambda$ is permanently irreducible.
\end{definition}

\begin{lemma}\label{never deleted}
A permanently irreducible word is permanently \RR-irreducible.
A permanent rule of \Store is never deleted.
A permanent rule is accepted by $\Rules[n+1]$ provided it is present in
\Store when the $n$-th Knuth--Bendix pass begins; it is then accepted
by $\Rules[m]$ for all $m > n$.
\end{lemma}
\begin{proof}
Let $u$ be permanently irreducible.
\RR-reduction of $u$ can only take place if, immediately after the
\RR-reduction, $u$ is \Store-reducible, conceivably as a result of
some rule being added to \Store during the \RR-reduction.
But this is impossible by hypothesis.

A rule $(\lambda,\rho)$ is deleted only as a result of being the input
to the minimization procedure.
By Lemma~\ref{maintain congruence}, there would have to be a Thue path
from $\lambda$ to $\rho$ with largest node $\lambda$. The first
elementary reduction must therefore be rightward (see
Definition~\ref{Thue path}) $\lambda\rightarrow_{\Storen t} \mu$.
We are assuming that $(\lambda,\rho)$ is a permanent rule of \Store.
Since every 
proper subword of $\lambda$ is permanently irreducible, it is
permanently \RR-irreducible, as we have just seen.
So this first
elementary reduction must be associated to a rule $(\lambda,\mu)$.

Either $\mu = \rho$, in which case the rule $(\lambda,\rho)$ has not
been deleted, or else,
when $(\lambda,\rho)$ was input to the minimization
routine, $\rho$ was \RR-reducible. However, it is permanently
\RR-irreducible which is a contradiction.

It follows that if $(\lambda,\rho)$ is present in \Store
at the start of the $n$-th
Knuth--Bendix pass, it will be sewn into \WDiff at some point during the $n$-th
Knuth-Bendix pass and accepted by $\Rules[n+1]$.
Since $(\lambda,\rho)$ is a permanent rule, it will subsequently
remain in \Store and will be presented for minimization during each pass.
The same rule will be output and used to mark
states and arrows of \WDiff as \needed.
Therefore, $(\lambda,\rho)$ is accepted by $\Rules[m]$ for each $m\geq n$. 
\end{proof}

\begin{lemma}\label{eventually permanent}
Let $u$ be a fixed word. Then there is a $t_0$ depending on $u$, such that,
for all $t \ge t_0$, each elementary \Storen t-reduction of $u$ is associated
to a permanent rule.
If all proper subwords of $u$ are permanently irreducible, then, for $t\ge t_0$,
there is at most one elementary reduction of $u$, and this is associated
to a permanent rule $(u,w)$. 
\end{lemma}
\begin{proof}
There are only finitely many subwords of $u$. So we need only prove
that, given any word $v$, there is a $t_0$ such that for all $t\ge t_0$,
each rule in \Storen t with left-hand side $v$ is permanent.
If there is a proper subword of $v$ which is not permanently irreducible,
then at some time $s_0$ it becomes \Storen{s_0}-reducible.
By Lemma~\ref{stays S-reducible}, it is \Storen s-reducible for $s\ge s_0$.
By Lemma~\ref{stays RR reducible}, it becomes \RR-reducible at the beginning of
the next Knuth--Bendix pass after $s_0$. During this pass all rules
with left-hand side $v$ will be deleted. Also, since this proper subword
of $v$ is now permanently \RR-reducible, no rule with left-hand side
equal to $v$ will ever be inserted subsequently. In this case, the
result claimed about $v$ is vacuously true.

So we assume that each proper subword of $v$ is permanently irreducible,
and that $v$ itself is \Store-reducible at some time $t$. A rule $(v,w)$ will 
be permanent if
$w$ is permanently irreducible. Otherwise it will disappear as a result
of minimization and, by Lemma~\ref{no repeats}, never reappear.
There cannot be two permanent rules $(v,w_1)$ and $(v,w_2)$ with $w_1 > w_2$.
For critical pair analysis would produce a new rule $(w_1,w_2)$ during
the next Knuth--Bendix pass, and so $w_1$ would not be permanently
irreducible.
\end{proof}

\begin{figure}[htbp]
\begin{center}
{\fontsize{10}{12pt}\selectfont
\begin{picture}(200,220)(25,-130)
\put(120,70){$h_0$}
\put(50,50){\framebox(150,5)}
\put(90,50){\framebox(50,5)}
\put(110,50){\framebox(50,5)}
\multiput(90,51)(0,2){2}{\line(1,0){50}}
\multiput(110,52)(0,2){2}{\line(1,0){50}}
\put(110,60){$\lambda_1$}
\put(135,60){$\lambda_2$}
\put(95,40){$u_1^\prime$}
\put(120,40){$u_2^\prime$}
\put(145,40){$u_3^\prime$}

\put(-40,0){\framebox(150,5){}}
\put(0,0){\framebox(50,5){}}
\put(50,0){\framebox(20,5){}}
\multiput(0,1)(0,2){2}{\line(1,0){50}}
\multiput(50,2)(0,2){2}{\line(1,0){20}}
\put(20,10){$\rho_1$}
\put(55,-10){$u_3^\prime$}
\put(110,35){\vector(-4,-1){70}}

\put(140,0){\framebox(150,5){}}
\put(180,0){\framebox(20,5){}}
\put(200,0){\framebox(50,5){}}
\multiput(180,1)(0,2){2}{\line(1,0){20}}
\multiput(200,2)(0,2){2}{\line(1,0){50}}
\put(185,-10){$u_1^\prime$}
\put(220,10){$\rho_2$}
\put(140,35){\vector(4,-1){70}}


\put(-40,-70){\framebox(150,5){}}
\put(0,-70){\framebox(70,5){}}
\multiput(0,-69)(0,1){4}{\line(1,0){70}}
\put(30,-60){$z_1$}
\put(140,-70){\framebox(150,5){}}
\put(180,-70){\framebox(70,5){}}
\multiput(180,-69)(0,1){4}{\line(1,0){70}}
\put(210,-60){$z_2$}

\put(30,-45){$\vee$}
\put(33,-45){\line(0,1){35}}
\put(37,-45){$\star$}
\put(7,-45){$\Storen {t^{\prime\prime}}$}

\put(210,-45){$\vee$}
\put(213,-45){\line(0,1){35}}
\put(217,-45){$\star$}
\put(187,-45){$\Storen {t^{\prime\prime}}$}

\qbezier(30,-80)(125,-160)(210,-80)
\put(35,-80){\line(-1,0){5}}
\put(30,-85){\line(0,1){5}}
\put(205,-80){\line(1,0){5}}
\put(210,-85){\line(0,1){5}}
\put(203,-78){$\star$}
\put(212,-90){$\Storen {t^{\prime\prime}}$}
\end{picture}
}
\end{center}
\caption{\sf When the leftward and rightward reductions from $h_0$ are obtained
from rules $(\lambda_1,\rho_1)$ and $(\lambda_2,\rho_2)$ having overlapping
left-hand sides, this diagram shows the time-$t^{\prime\prime}$ Thue
paths that exist after the resulting critical pair analysis.}
\label{non-disjoint left-hand sides}
\end{figure}
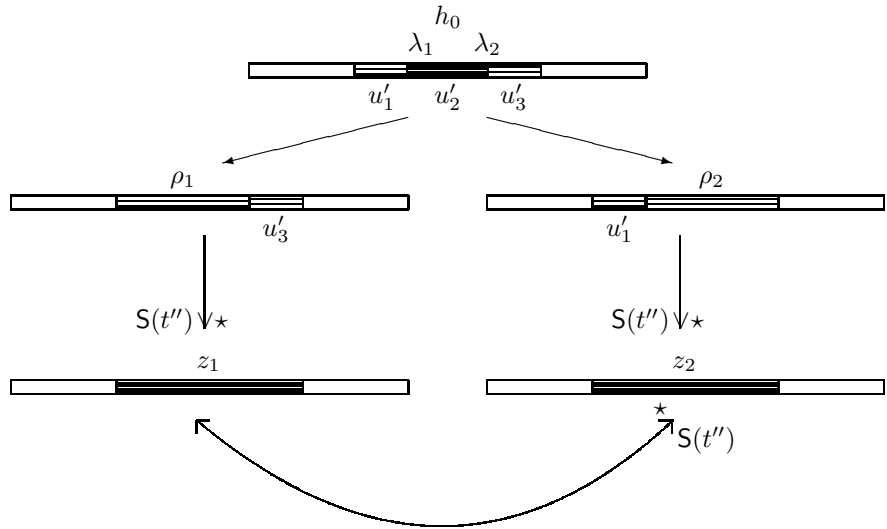

\begin{theorem}\label{confluence}
Let $u$ be a fixed word in $\Astar$ and let $v$ be the smallest
element in its Thue congruence class. Then, for large enough times,
there is a chain of elementary reductions from $u$ to $v$ each associated
to a permanent rule.
After enough time has
elapsed, \RR-reduction of $u$ always gives $v$.
(Recall that $v$ is the \shl representative of $\overline u$.)
\end{theorem}
\begin{proof}
We start by proving the first assertion.
By hypothesis, we have, for each time $t$, a time-$t$ Thue path
$p_t$ from $u$ to $v$, and we can suppose that $p_t$ contains no repeated
nodes by cutting out part of the path if necessary.
The only reason why we couldn't take $p_{t+1}$ to be $p_t$ is if some
rule $(\lambda,\rho)$, used along the Thue path $p_t$, is deleted at
time $t$. By Lemma~\ref{maintain congruence} we can, however, assume that
each node of $p_{t+1}$ is either already a node of $p_t$ or is smaller
than some node of $p_t$.

Let $h_0$ be the largest node on $p_0$, and suppose that we have already
proved the theorem for all pairs $u$ and $v$ which are connected by a Thue
path with largest node smaller than $h_0$. By induction on $t$, using
\thref{maintain congruence}, we can assume
that $h_0$ is the largest node on $p_t$ for all time $t$. If $v=h_0$
then since $v$ is the smallest element in its congruence class, there
are no elementary reductions starting from $v$, and we must have $u=v$
in this case.

By Lemma~\ref{eventually permanent}, we may assume that $t_0$ has been
chosen with the property that, for all words $w\le h_0$ and for all
$t\ge t_0$, all elementary \Storen t-reductions of $w$ are associated
to permanent rules which are accepted by $\Rules[n]$ provided $n$ is
sufficiently large.

Let $h_0 = \mu_t\alpha_t \nu_t\to_{\Storen t} \mu_t\beta_t \nu_t$ be the
rightward elementary reduction of $h_0$ at time $t$.
Our construction of $p_{t+1}$ from $p_t$, as in \thref{maintain
congruence}, makes $\alpha_{t+1}$ a subword of $\alpha_t$.
The construction also ensures that,
if $\alpha_{t+1} = \alpha_t$, then $\beta_{t+1} \le \beta_t$.
The rule $(\alpha_t,\beta_t)$ is therefore independent of $t$ for large values
of $t$.
Then $(\alpha_t,\beta_t)$
is permanent and $\alpha_t$ is \RR-reducible for large enough $t$.
If $u\neq h_0$, the same argument applies to
the unique elementary leftward reduction with source $h_0$ at time $t$.

If $h_0 = u$, let $u \rightarrow_{\Storen t} w$ be the first rightward 
elementary reduction for large
values of $t$. By our induction hypothesis, there is a Thue path of
elementary reductions from $w$ to $v$, each associated to a permanent
rule, and with no node larger than $w$, and so we have the required Thue
path from $u$ to $v$. 

Suppose now that $h_0\neq u$, so that we get two permanent rules,
associated to the leftward and rightward elementary reductions of
$h_0$.
If the two elementary reductions are identical, that is, if the two
permanent rules are equal and if their left-hand sides occur in the
same position in $h_0$, then $p_t$ contains a repeated node which we are
assuming not to be the case. So the two elementary reductions occur in
different positions in $h_0$. Now choose $t$ to be large enough so that the two
rules concerned have already been compared in a critical pair analysis
in Step~\ref{compare for overlaps} during some previous Knuth--Bendix
pass.

If these two rules have left-hand sides which are
disjoint subwords of $h_0$, then we can interchange their order so as to
obtain a Thue path from $u$ to $v$ where all nodes are strictly smaller than
$h_0$---see Figure~\ref{disjoint left-hand sides}. The first assertion of
the theorem then follows by the induction hypotheses in this
particular case.

\begin{figure}[htbp]
\begin{center}
{\fontsize{10}{12pt}\selectfont
\begin{picture}(200,150)(25,-50)
\put(50,50){\framebox(150,5){}}
\put(80,50){\framebox(30,5){}}
\put(140,50){\framebox(30,5){}}
\multiput(80,51)(0,1){4}{\line(1,0){30}}
\multiput(140,51)(0,1){4}{\line(1,0){30}}
\put(90,60){$\lambda_1$}
\put(150,60){$\lambda_2$}
\put(120,70){$h_0$}
\put(80,40){\vector(-4,-1){70}}
\put(-50,5){\framebox(150,5){}}
\put(-20,5){\framebox(30,5){}}
\put(40,5){\framebox(30,5){}}
\multiput(-20,6)(0,1){4}{\line(1,0){30}}
\multiput(40,6)(0,1){4}{\line(1,0){30}}
\put(-10,15){$\rho_1$}
\put(50,15){$\lambda_2$}
\put(170,40){\vector(4,-1){70}}
\put(150,5){\framebox(150,5){}}
\put(180,5){\framebox(30,5){}}
\put(240,5){\framebox(30,5){}}
\multiput(180,6)(0,1){4}{\line(1,0){30}}
\multiput(240,6)(0,1){4}{\line(1,0){30}}
\put(190,15){$\lambda_1$}
\put(250,15){$\rho_2$}
\put(50,-40){\framebox(150,5){}}
\put(80,-40){\framebox(30,5){}}
\put(140,-40){\framebox(30,5){}}
\multiput(80,-39)(0,1){4}{\line(1,0){30}}
\multiput(140,-39)(0,1){4}{\line(1,0){30}}
\put(70,-5){\vector(4,-1){70}}
\put(180,-5){\vector(-4,-1){70}}
\put(90,-30){$\rho_1$}
\put(150,-30){$\rho_2$}
\put(40,30){\dashbox{5}(170,50){}}
\end{picture}
}
\end{center}
\caption{\sf Removing the node $h_0$ when the leftward and rightward reductions
are obtained from rules having disjoint left-hand sides.}
\label{disjoint left-hand sides}
\end{figure}
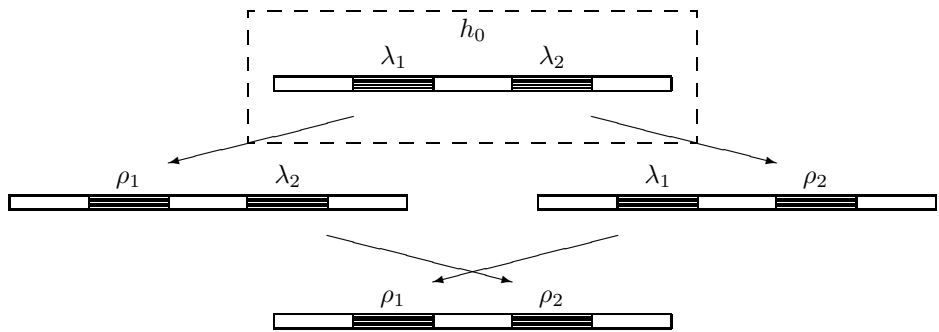

If the two left-hand sides do not correspond to disjoint subwords of
$h_0$ then, by assumption, there is some time
$t^\prime < t$, such that a critical pair
$(u^\prime,v^\prime,w^\prime)$ was considered.
Here
$u^\prime \rightarrow_{\Storen {t^\prime}} v^\prime$ and
$u^\prime \rightarrow_{\Storen {t^\prime}} w^\prime$ are elementary
$\Storen {t^\prime}$-reductions given by the two rules, and $u^\prime$ is a
subword of $h_0$. After the critical pair analysis, at time
$t^{\prime\prime}\leq t$, the Thue paths illustrated in
Figure~\ref{non-disjoint left-hand sides} are possible.
As a consequence of
\thref{maintain congruence}, it is straightforward to see that for all times
$s\geq t^{\prime\prime}$, $v^\prime$ and $w^\prime$ can be connected by a
time-$s$ Thue path
in which all nodes are no larger than the largest of $v^\prime$ and $w^\prime$.
In particular, this applies at time $t$ so that the targets of the two
elementary
\Storen t-reductions from $h_0$ can be connected by a time-$t$ Thue path in
which all nodes are strictly smaller than $h_0$.
This completes the inductive proof of the first assertion of the
theorem.

We have arranged that $t$ is large enough so that,
for all $w\le u$, all elementary \Storen t-reductions of $w$ are
associated to permanent rules, and such a $w$ can be permanently \RR-reduced
to the least element in its Thue congruence class. It follows that such a $w$ 
is \RR-irreducible if and only if it is minimal in its Thue class. In 
particular \RR-reduction of $u$ must give $v$.
\end{proof}

\begin{corollary}\label{rule corollary}
(i) The set of permanent rules in \RR is confluent. (ii)
The set of
such rules is equal to $\Permanent = \bigcap_{t}\bigcup_{s\ge t}$\Storen s.
(iii) A word $u$ is smallest in its Thue congruence class if
and only if it is permanently irreducible and this is equivalent to
being in \shl normal form.
(iv) Each permanent rule is a $U$-minimal rule and each $U$-minimal
rule is accepted by $\Rules[n]$ for $n$ sufficiently large.
\end{corollary}
\begin{proof}
The first and third statements are obvious from Theorem~\ref{confluence}.
For the second statement, each permanent rule is contained in \Permanent by
Lemma~\ref{never deleted}. Conversely, if we have a rule $r$ in \Store which is
not permanent, then for all
sufficiently large times $s$ either its right-hand side or a proper subword
of its left-hand side is \Storen s-reducible. Theorem~\ref{confluence} ensures
that this reducible word is \RR-reducible for all sufficiently large
times $s$. Therefore $r$ will be minimized and deleted from
\Store. Hence from Lemma~\ref{no repeats} we see that $r$ is not contained in
\Permanent.

To prove the fourth statement, suppose $(\lambda,\rho)$ is
$U$-minimal.
By \thref{confluence}, a Thue path from $\lambda$ to $\rho$
will eventually be generated by our Knuth--Bendix procedure and each
elementary reduction in the path will be rightward and
associated to a permanent rule.
The first elementary reduction must have the form $(\lambda,\rho')$,
because each proper subword of $\lambda$ is permanently irreducible.
But then $\rho'=\rho$, for otherwise $\rho'>\rho$ and
\thref{confluence} applies to show that $\rho'$ is not permanently
irreducible. But then $(\lambda,\rho')$ would not have been a
permanent rule. Therefore $(\lambda,\rho)$ is a permanent rule.

Conversely, suppose that $(\lambda,\rho)$ is a permanent rule.
This means that $\rho$ and every proper subword of $\lambda$ is
permanently irreducible.
By \thref{confluence}, this mens that $\rho$ and every proper subword
of $\lambda$ are in \shl normal form.
It follows that $(\lambda,\rho)$ is $U$-minimal.
\end{proof}

The next result is the main theorem of this paper.
\begin{theorem}\label{main}[Theorem]
Let $G$ be a group with a given finite presentation
and a given ordering of the generators and their inverses.
Suppose that the set of $U$-minimal rules is regular (for example if 
$(G,A)$ is \shl-automatic).
Then the procedure given in \ref{main loop} will stabilize at some $n_0$ with
$\Rules[n+1] = \Rules[n]$ if $n\ge n_0$.
\Permanent (defined in \thref{rule corollary})
is then the language of a certain
two-variable finite state automaton and the automaton can be explicitly
constructed. (Unfortunately we do not have a method of knowing when or
whether we have reached $n_{0}$.)
\end{theorem}
\begin{proof}
By hypothesis there is a two-variable automaton accepting the set of all
$U$-minimal rules. By welding, we obtain a two-variable rule automaton $M$.
By amalgamating states, we
may assume that each state of $M$ corresponds to a different
word-difference.

Given any arrow in $M$, there is a $U$-minimal rule $(\lambda,\rho)$
which is accepted by $M$ and which uses that arrow.
By \thref{rule corollary}. $(\lambda,\rho)$ is a permanent rule which
is eventually generated by our Knuth--Bendix procedure. 
By \thref{never deleted}, such a rule is never deleted.
Since there are only a finite number of arrows in $M$, we see that,
for large enough $n$, each $(\lambda,\rho)$ in this finite set of
rules may be traced out in $\Rules[n]$.
We record the states and arrows reached as being required by this finite set
of rules.

We may also assume that the
states in $\Rules[n]$ which have been recorded as just explained,
are all associated to different word-differences. To see this,
first note that any equality of word-differences between different states
is eventually discovered according to \thref{confluence}.
Then, as in \ref{WDiffdetails}, the corresponding states are
amalgamated. It follows that, for $n$ large enough, there is a copy of
$M$ inside $\Rules[n]$.

Subsequently, arrows and states lying outside $M$
will not be used in \RR-reduction. They will not be
marked as \needed and will be deleted. It follows that
$\Rules[n] =M$ for $n$ sufficiently large.

Finally, knowing $M$, we can easily change it to a finite state
automaton accepting exactly the minimal rules---this involves making
sure that if $(u,v)$ is accepted, then $u>v$, $v$ is irreducible and
every proper subword of $u$ is irreducible.
\end{proof}

%% file: Sec7.tex
\section{Fast reduction}
\label{Fast reduction}[Section]

In this section, we show how to rapidly reduce an arbitrary word,
using the rules in \SetOfRules\Rules together with the rules in \Store.
We assume the properties made explicit in \ref{properties of Rules}.
The time taken to carry out the first reduction
is bounded by a small constant times the length of the word.
This efficiency is possible because of the use of finite state
automata to do the reduction.

\subsection{Rules for which no prefix or suffix is a rule.}
At the moment, it is possible for an element $(u,v)^+$ of $\SetOfRules\Rules$
to have a prefix or suffix which is also a rule. This is undesirable because
it makes the computations we will have to do bigger and longer without any
compensating gain.

Recall that the automaton recognizing $\SetOfRules\Rules$
is the product of $\Rules$ with $\SLtwo$, the initial state being the product
of initial states and the set of final states being any product of final states.
By \ref{properties of Rules}, there is only one initial and one final
state of \Rules; these are equal and the state is denoted by $s_0$.

We remove from $\Rules$ any arrow labelled $(x,x)$ from the initial state to
itself. We then form the product automaton, as described above, with two
restrictions. Firstly, we omit any arrow whose source is a product of final
states. Secondly, we omit the state
with first component equal to $s_0$, the initial state of
$\Rules$, and second component equal to state $3$ of $\SLtwo$
(see Figure~\ref{SL2}) and any arrow
whose source or target is this omitted state.
We call the resulting automaton $\Rules^\prime$.

\begin{lemma}\label{arrows removed 1} The language accepted by
$\Rules^\prime$ is the set of labels
of accepted paths in the product automaton, starting from the product
of initial states and ending at a product of final states, such that the only
states along the path with first component equal to $s_0$ are at the
beginning and end of the path.
\end{lemma}
\begin{proof} First consider an accepted path $\alpha$ in $\Rules'$.
The only arrows in $\Rules'$ with source having first component $s_0$
are those with source the product of initial states.
In $\SLtwo$ it is not possible to
return to the initial state. It follows that $\alpha$ has the required form.

Conversely any such path in the product automaton also lies in $\Rules'$
because it avoids all omitted arrows.
\end{proof}

\begin{lemma}\label{arrows removed 2} The language accepted by
$\Rules^\prime$ is the subset of $\SetOfRules\Rules$ which has no proper
suffix or proper prefix in $\SetOfRules\Rules$.
\end{lemma}
\begin{proof}
If $\alpha$ is an accepted path in $\Rules'$, then it is clearly in
$\SetOfRules\Rules$. Moreover if it had a proper suffix or proper prefix
which was in $\SetOfRules\Rules$, there would be a state in the middle
of $\alpha$ with first component $s_0$. We have seen that this is impossible
in Lemma~\ref{arrows removed 1}.

Conversely, we must show that if $\alpha$ is an accepted path in
the product automaton such that no proper prefix and no proper suffix of
$\alpha$ would be accepted by the product automaton, then no state
met by $\alpha$, apart from its two ends, has $s_0$ as a first component.
Let $\alpha = ((s_0,1),(u_1,v_1), q_1, \ldots, (u_n,v_n), q_n)$,

First suppose $u_1 < v_1$. Since $\alpha$ is accepted by $\SLtwo$,
$|u| > |v|$ and we must have $v_n = \$$.
Let $r<n$ be chosen as large as possible so that the first component
of $q_r$ is $s_0$.
Then $(u_{r+1},v_{r+1})\ldots (u_n,v_n)$ will be accepted by $\Rules$
and will be accepted by $\SLtwo$ because $v_n= \$$.
Since this cannot be a proper suffix of $\alpha$ by assumption, we must have
$r=0$. Hence $q_i$ has a first component equal to $s_0$ if and only if
$i=0 $ or $i=n$.

Next note that we cannot have $u_1=v_1$. This is because there is no
arrow labelled $(u_1,u_1)$ in $\SLtwo$ with source the initial state,
so $\alpha$ would not be accepted by the product automaton.

Now suppose that $u_1>v_1$ and let $r>0$ be chosen as small as
possible so that the first component of $q_r$ is $s_0$.
Since $u_1>v_1$, the second component of $q_r$ will be a final state
(see Figure~\ref{SL2}). Since $\alpha$ has no accepted proper prefix, we
must have $r=n$. Hence $q_i$ has a first component equal to $s_0$ if
and only if $i=0 $ or $i=n$.

So we have proved the required result for each of the three
possibilities.
\end{proof}

Reduction with respect to $\SetOfRules\Rules$
is done in a number of steps. First we find the shortest reducible
prefix of $w$, if this exists. Then we find the shortest suffix of that which
is reducible. This is a left-hand side of some rule in $\SetOfRules\Rules$.
Then we find the corresponding right-hand side and substitute this for the
left-hand side which we have found in $w$. This reduces $w$ in the
\shl-order. We then repeat the operation until we obtain an
irreducible word. The process is explained in more detail in
\ref{history stack}.

Our first objective
is to find the shortest reducible prefix of $w$, if this exists.
To achieve this, we must determine whether $w$ contains a
subword which is the left-hand side of rule belonging to $\SetOfRules\Rules$.

Let $\Rules^{\prime\prime}$ be the automaton obtained
from $\Rules^\prime$ (see Lemmas~\ref{arrows removed 1} and
\ref{arrows removed 2}) by adding
arrows labelled $(x,x)$ from the initial state to the initial state.

We construct an FSA $\NRed\Rules$ in one variable
by replacing each label of the form $(x,y)$
on an arrow of $\Rules^{\prime\prime}$ by $x$.
Here $x\in A$ and $y\in A^+$.
The name of the automaton $\NRed\Rules$
refers to the fact that the automaton accepts reducible
words, and does so non-deterministically.
We obtain an FSA with no $\epsilon$-arrows.
However there may be many arrows labelled $x$ with a given source.
Let $\LHS\Rules$ be the regular language of left-hand sides of rules in
$\SetOfRules\Rules$ such that no proper prefix or proper suffix of the rule
is itself a rule.

\begin{lemma}
$\Astar.\LHS\Rules= L(\NRed\Rules)$.
\end{lemma}
\begin{proof}
Because of the extra arrows labelled $(x,x)$ from initial state to
initial state, inserted into $\Rules^{\prime\prime}$, the inclusion
$\Astar.\LHS\Rules\subset L(\NRed\Rules)$ is clear.

Conversely, if $u$ is accepted by $\NRed\Rules$, there is a corresponding pair
$(u,v)$ accepted by $\Rules^{\prime\prime}$. We find a maximal common prefix
$p$ of $u$ and $v$, so that $u=pu'$ and $v=pv'$.
$\Rules^{\prime\prime}$ remains in the initial state while reading $(p,p)$.
Since the initial state of $\SLtwo$ is not a final state, $(u',v')$
must be non-empty. Since there is no way of
returning to the initial state of $\SLtwo$, once $\Rules^{\prime\prime}$ starts
reading $(u',v')$, it can never return to the initial state, and therefore
$(u',v')$ must be accepted by $\Rules^\prime$. Therefore $u'\in \LHS\Rules$,
as claimed.
\end{proof}

\subsection{The automaton $P$.} \label{P}
To find the shortest reducible prefix of a given word $w$
we could feed $w$ into the FSA $\NRed\Rules$. However, reading a word
with a non-deterministic automaton is very time-consuming, as all
possible alternative paths need to be followed.

For this reason, it may at first sight seem sensible to
determinize the automaton. However, determinizing a non-deterministic
automaton potentially leads to an exponential increase in size.
The states of the determinized automaton are subsets of the non-deterministic
automaton, and there are potentially $2^n$ of them if there were $n$
states in the non-deterministic automaton.

For this reason, we use a {\it lazy
state-evaluation} form of the subset construction. The lazy evaluation strategy
(common in compiler design---see for example \cite{AhoSethiUllman})
calculates the
arrows and subsets as and when they are needed, so that a gradually increasing
portion $P(\Rules)$ of a determinized version $\Red\Rules$ of
$\NRed\Rules$ is all that
exists at any particular time.

Lazy evaluation is not automatically an advantage. For example, if in the
end one has to construct virtually the whole determinized automaton
$\Red\Rules$ in any
case, then nothing would be lost by doing this immediately. In our special
situation,
lazy evaluation \textit{is} an advantage for two reasons. First,
during a single pass of the Knuth--Bendix process
(see \ref{Knuth--Bendix pass}), only a
comparatively small part of the determinized one-variable automaton
$\Red\Rules$ needs to be constructed.
In practice, this phenomenon is particularly
marked in the early stages of the computation, when the automata are
far from being the ``right'' ones. Second,
this approach gives us the
opportunity to abort a pass of Knuth--Bendix, recalculate on the basis
of what has been discovered so far in this pass, and then restart the
pass. If an abort seems advantageous early in the pass, very little
work will have been done in making the structure of a determinized
version of $\Red\Rules$ explicit.

At the start of a Knuth--Bendix pass
we let $P(\Rules)$ be the one-variable automaton containing only one 
state and no arrows. The state is an initial state of $P(\Rules)$ 
which is a singleton set whose only element is
the ordered pair of initial states of $\Rules$
and $\SLtwo$. At a subsequent time during the pass, $P(\Rules)$ may have increased,
but it will always be a portion of $\Red\Rules$.
Each state of $P(\Rules)$ is a set of pairs $(s,t)$, where $s$ is a 
state of \Rules and $t$ is a state of \SLtwo.

The transition with source $s$, a state in $P(\Rules)$, and label $x\in A$ may
or may not already be defined.
If it is defined, we denote by $\mu(s,x)$ the target of this arrow. 

Suppose now that we wish to find the shortest prefix of the word
$w=x_1\cdots x_n\in A\uast$ which is $\SetOfRules\Rules$-reducible.
Suppose that $s_0,s_1,\ldots,s_k$ are states of $P(\Rules)$, where
$0\leq k\leq n-1$,
that $s_0$ is the start state of $P(\Rules)$,
and that,
for each $i$ with $1\leq i\leq k$, the arrow with source $s_{i-1}$ and 
label $x_{i}$ has been constructed, with target
$\mu(s_{i-1},x_i) = s_i$.
Suppose that
the target of the arrow with source $s_{k}$ and label $x_{k+1}$
has not yet been defined.

The conventional subset construction applied to the state $s_k$ of
$P(\Rules)$ under the alphabet symbol $x_{k+1}$ yields a set, which
we denote by $\mu_1(s_k,x_{k+1})$.
This is how $\mu_1(s_k,x_{k+1})$ is defined.
For each $(s',t')\in s_k$, we look
for all arrows in $\NRed\Rules$
labelled $x_{k+1}$ with source $(s',t')$.
If $(s,t)$ is the target of such an arrow, then $(s,t)$ is an element
of $\mu_1(s_k,x_{k+1})$.
Note that this subset is always non-empty, because the initial
state of $\NRed\Rules$ is an element of each $s_i$.

In the standard determinization procedure
one would now look to see whether there is already a state $s_{k+1}$ of
$P(\Rules)$ which is equal to $\mu_1(s_k,x_{k+1})$.
If not, one would create such a state $s_{k+1}$.
One would then insert an
arrow labelled $x_{i+1}$ from $s_k$ to $s_{k+1}$, if there wasn't
already such an arrow.
A new state is defined to be a final state of $P(\Rules)$ if and
only if the subset contains a final state of $\NRed\Rules$.
Of course, one does not need to determine the subset $\mu_1(s_k,x_{k+1})$
if there is already an arrow in $P(\Rules)$ labelled $x_{k+1}$ with
source $s_k$, because in that case the subset is already computed and stored.

In our procedure we improve on the procedure just
described. The point is that
$\mu_1(s_k,x_{k+1})$ may contain pairs which are
not needed and can be removed. From a practical point of view
this has the advantage of saving space and reducing the amount of computation
involved when calculating subsequent arrows. Specifically, we remove a
pair $(p,q^\prime)$ from $\mu_1(s_k,x_{k+1})$ if
$q^\prime$ is state $3$ of $\SLtwo$ (see Figure~\ref{SL2}) and
$\mu_1(s_k,x_{k+1})$ also contains the pair $(p,q)$
where $q$ is state $2$ of \SLtwo (same $p$ as in $(p,q^\prime)$)
Removing all such pairs
$(p,q^\prime)$ yields the set $\mu_P(s_k,x_{k+1})$ and we add the
corresponding arrow and state to $P(\Rules)$,
creating a new state if necessary. We make the state a final
state if the subset contains a final state of $\NRed\Rules$. The
validity of this modification follows from
Theorem~\ref{modified determinization}, and we see that
some prefix of $w$ arrives at a final state of $P(\Rules)$ if and only if
$w$ is $\SetOfRules\Rules$-reducible.

When finding the corresponding left-hand side of a rule inside $w$,
we need never compute beyond a final state of $P(\Rules)$.
As a space-saving and time-saving measure our
implementation therefore replaces each final state of $P(\Rules)$, as soon as
it is found, by the empty set of states. As remarked above, the standard
determinization
of $\NRed\Rules$ never produces an empty set of states, so there is no
possibility of confusion.

Reading $w$ can be quite slow if many states need to be added to $P(\Rules)$
while it is being read. However, reading $w$ is fast when no states
need to be built. In practice, fairly soon after a Knuth--Bendix pass
starts, reading becomes rapid, that is, linear with a very small constant.

\subsection{Finding the left-hand side in a word.}
\label{finding the left-hand side}
We retain the hypotheses of Section~\ref {Fast reduction}.
Namely, we have a two-variable automaton $\Rules$ satisfying the
conditions of Paragraph~\ref{properties of Rules}. We are given a word
$w=x_1\cdots x_n$, and we wish to reduce it. In the previous section
we showed how to find the minimal reducible prefix $w'=x_1\cdots x_m$
of $w$ with respect
to the rules implicitly specified by $\Rules$. We now wish to find the
minimal suffix of $w'$ which is a left-hand side of some rule in
$\SetOfRules\Rules$. The procedure is quite similar to that of the
previous section.

We will now give the basic construction.
However, the details will later need to be modified so as to achieve
greater computational efficiency in finding the associated right-hand side,
if this is necessary.
Our reason for including the simpler version is to lead the reader
more gently and with more understanding to the actual more complex version.

\begin{figure}[htbp]
\begin{center}
{\fontsize{10}{12pt}\selectfont
\begin{picture}(250,160)(0,-50)
\put(160,80){\vector(0,-1){10}}
\put(150,50){\grid(20,20)(20,20)}
\put(154,62){$2\,4$}
\put(158,52){$5$}
\put(170,60){\vector(1,0){50}}
\put(183,65){$(x,y)$}
\put(182,52){$x>y$}
\put(220,50){\grid(20,20)(20,20)}
\put(222,52){\grid(16,16)(16,16)}
\put(224,58){$1\,2$}
\qbezier(225,70)(260,90)(240,55)
\put(239,74.5){$>$}
\put(235,93){$(x,y)$}
\put(235,83){$x>y$}
\put(220,-20){\grid(20,20)(20,20)}
\put(226,50){\vector(0,-1){50}}
\put(198,30){$(x,y)$}
\put(198,20){$x\leq y$}
\put(234,0){\vector(0,1){50}}
\put(240,30){$(x,y)$}
\put(240,20){$x>y$}
\put(227,-12){$2$}
\qbezier(225,-20)(260,-40)(240,-5)
\put(238,-29.5){$>$}
\put(233,-39){$(x,y)$}
\put(233,-49){$x\leq y$}
\qbezier(160,50)(160,-10)(220,-10)
\put(213,-12.5){$>$}
\put(137,15){$(x,y)$}
\put(140,3){$x\leq y$}

\put(150,60){\vector(-1,0){50}}
\put(76,48){\grid(24,24)(24,24)}
\put(78,50){\grid(20,20)(20,20)}
\put(82,62){$1\,2$}
\put(82,52){$3\,4$}
\put(78,-20){\grid(20,20)(20,20)}
\put(82,-12){$2\,3$}
\qbezier(83,-20)(118,-40)(98,-5)
\put(96,-29.5){$>$}
\put(91,-39){$(x,x)$}
\put(2,48){\grid(24,24)(24,24)}
\put(4,50){\grid(20,20)(20,20)}
\put(76,60){\vector(-1,0){50}}\
\put(8,62){$1\,2$}
\put(12,52){$3$}
\qbezier(2,53)(-18,92)(21,72)
\put(-2,76){$<$}
\put(-8,96){$(x,y)$}
\put(-8,86){$x\neq y$}
\qbezier(26,53)(81,53)(81,0)
\put(25.5,50.5){$<$}
\qbezier(14,48)(14,-10)(78,-10)
\put(71,-12.5){$>$}
\put(95,48){\vector(0,-1){48}}
\put(100,22){$(x,x)$}
\put(113,65){$(x,\$)$}
\put(39,65){$(x,\$)$}
\put(29,78){$(x,y),x\neq y$}
\put(40,33){$(x,y)$}
\put(45,23){$x\neq y$}
\put(0,0){$(x,x)$}
\end{picture}
}
\end{center}
\caption{\sf This PDFA arises by applying the accessible subset construction to
$\Rev\SLtwo$ in the case where the base alphabet has more than one element.
Each state is a subset of the state set of $\Rev\SLtwo$ and final states have
a double border. This PDFA, when reading a
pair $(u,v)$ from right to left, keeps track of whether $u$ is longer than $v$
or not, which it discovers immediately since padding symbols if any must occur
at the right-hand end of $v$. Note that this automaton is minimized.}
\label{RevSLsubsets}
\end{figure}
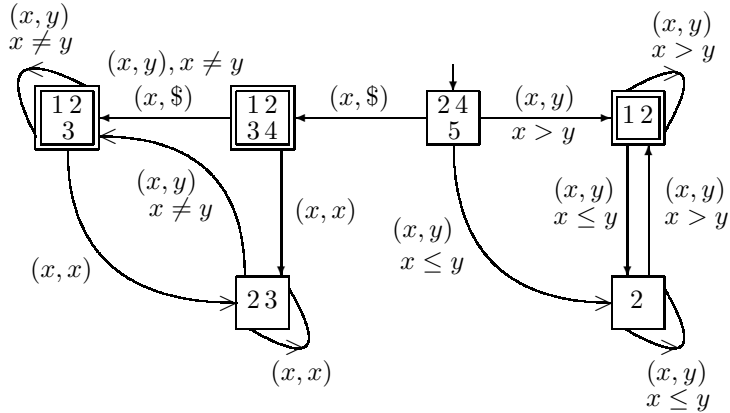

We form the two-variable
automaton $\Rev\Rules$, which we combine with $\Rev\SLtwo$.
The first automaton is, by hypothesis, partially deterministic.
If we determinize the second automaton, we obtain another PDFA.
Figure~\ref{RevSLsubsets} shows the determinization of $\Rev\SLtwo$,
where the subsets of states of $\SLtwo$ are explicitly recorded.

We take the product of the two automata
$\Rev\Rules$ and $\Rev\SLtwo$. A new state is a pair of old states.
An arrow is a pair of arrows with the same label $(x,y)$.
The initial state in the product is the unique pair of initial states.
A final state in the product is a pair of final states.

To form the one-variable non-deterministic automaton $\NRevLHS\Rules$
without $\epsilon$-arrows, we use the same states and arrows as in the
product automaton, but replace each label
of the form $(x,y)$ in the product automaton by
the label $x$.
The deterministic one-variable automaton $\DRevLHS\Rules$
can then be constructed using the subset construction.

As we have already warned the reader, we
use not the construction just described, but
a related construction which we describe below.
The point of what we do may not become fully apparent until we get to
\ref{finding the right-hand side}.

\subsection{Reversing the rules.}
\label{structure} We first describe a two-variable PDFA $M$ which accepts
exactly the reverse of each rule $(\lambda,\rho)^+$ in $\SetOfRules\Rules$
such that no proper suffix and no proper prefix of $(\lambda,\rho)^+$ is in
$\SetOfRules\Rules$ (cf. Lemma~\ref{arrows removed 2}).
We assume that we have a two-variable automaton $\Rules$ satisfying the
conditions of Paragraph~\ref{properties of Rules}.

A state of $M$ is a triple $(s,i,j)$, where $s$ is a state of $\Rev\Rules$,
$i\in\{0,1,2\}$ and $j\in\{+,-\}$.
The intention is that in a state $(s,i,j)$, $i$ represents the number of
padded symbols occurring in any path of arrows from the initial
state of $M$ to $(s,i,j)$. By \ref{right-hand sides at most
two shorter}, the padded symbols must be of the form $(x,\$)$, where $x\in A$.
There are zero, one or two padded symbols in any rule, and, if
padded symbols appear, they are at the right-hand end of a rule. This means
that they are the first symbols read by $M$.
The $j$ component is
intended to represent whether an arrow is permitted with source $(s,i,j)$
and label a padded symbol.
We take $j=+$ if a padded symbol is permitted, and $j=-$ if a padded
symbol is not permitted.

$M$ has a unique initial state $(s_0,0,+)$ where $s_0$
is the unique initial state of $\Rev\Rules$. In addition,
$M$ has three final states $f_0=(s_0,0,-), f_1=(s_0,1,-)$ and $f_2=(s_0,2,-)$.
We do not allow states of $M$ of the form $(s_0,i,j)$, except for the initial
state and the three final states just mentioned. We will construct the arrows
of $M$ to ensure that
any path of arrows accepted by $M$ has first component equal to $s_0$ for
its initial state and its final state and for no other states. (Compare this
with Lemma~\ref{arrows removed 1}.)

The following conditions determine the arrows in $M$.
\begin{enumerate}
\item Each arrow of $M$ is labelled with some $(x,y)$, where
$x\in A$ and $y\in A^+$.
\item $(s,i,j)^{(x,\$)}$ is defined if and only if 1) $t=s^{(x,\$)}$ is defined
in $\Rev\Rules$, and 2a) $(s,i,j)=(s_0,0,+)$, the initial state, or
2b) $(i,j)=(1,+)$. In case 2a) the
target is $(t,1,+)$,
unless $t$ is the final state of $\Rev\Rules$, in which case
the target is $f_1 = (s_0,1,-)$. In case 2b), the target is $(t,2,-)$,
which may possibly be equal to $f_2$.
The final state $f_1$ arises in case 2a) when we have a rule $(x,\epsilon)$,
which means that the generator $x$ of our group represents the trivial element.
The final state $f_2$ arises in case 2b) when we have a rule
$(x_1x_2,\epsilon)$. This kind of rule arises when $x_1$ and $x_2$ are inverse
to each other, usually formal inverses.
\item For $i=0,1,2$, there are no arrows with source $f_i$.
\item Suppose $(s,i,j)$ is not a final state.
\label{non-padded M transitions}
Then $(s,i,j)^{(x,y)}$ with $x,y\in A$ is
defined if and only if
1) $t=s^{(x,y)}$ is defined in $\Rev\Rules$, and 2) if $t=s_0$ then
2a) $i=0$ and $x>y$ or 2b) $i>0$ and $x\neq y$.
We then have $(s,i,j)^{(x,y)} = (t,i,-)$.
This condition corresponds to the requirement that $(u,v)$ can only be a rule
if a) $u$ and $v$ have the same length and $u_1 > v_1$, where these are the
first letters of $u$ and $v$ respectively, or b)
if $u$ is longer than $v$ and $u_1\neq v_1$.
\end{enumerate}

\begin{lemma}
The language accepted by $M$ is the set of reversals of rules
$(\lambda,\rho)^+ \in \SetOfRules\Rules$ such that no proper suffix
and no proper prefix of
$(\lambda,\rho)^+$ is in $\SetOfRules\Rules$.
\end{lemma}
The proof of this lemma is much the same as the proofs of
Lemmas~\ref{arrows removed 1} and \ref{arrows removed 2}.
We therefore omit it.

Using the above description of $M$, we now describe how to obtain a
non-deterministic
one-variable automaton $\NRevLHS\Rules$ from $M$ in an analogous
manner to that
used to obtain $\NRed\Rules$ from $\Rules''$ in
Section~\ref{Fast reduction}.
$\NRevLHS\Rules$ accepts reversed left-hand sides
of rules in $\SetOfRules\Rules$ which do not have a proper prefix or a proper
suffix which is in $\SetOfRules\Rules$.
$\NRevLHS\Rules$ has the same set of states as $M$ and the same set of
arrows. However, the label $(x,y)$ with $x\in A$ and $y\in A^+$
of an arrow in $M$ is replaced by the label $x$ in
$\NRevLHS\Rules$
The two automata, $M$ and $\NRevLHS\Rules$, have the same initial state
and the same final
states. Hence $\NRevLHS\Rules$ accepts all reversed left-hand
sides $\lambda^R$ of rules $(\lambda,\rho)$
whose reversals $((\lambda,\rho)^+)^R$ are accepted by $M$.

\subsection{The automaton $Q$.}\label{Q}
The one-variable automaton $Q(\Rules)$ is
formed from $\NRevLHS\Rules$ by a modified subset construction, using lazy
evaluation. $Q(\Rules)$ is part of the one-variable
PDFA $\DRevLHS\Rules$, the determinization of $\NRevLHS\Rules$. As we shall
see, a word is accepted by $Q(\Rules)$ only if its reversal $\lambda$ is the
left-hand side of a rule in $\SetOfRules\Rules$ and no proper subword of
$\lambda$ has this property.

\begin{note}
In order to
construct states and arrows in $Q(\Rules)$, one only needs to have
access to $\Rev\Rules$, that is, neither $M$ nor $\NRevLHS\Rules$ has to be
explicitly constructed.
\end{note}

\subsection{The algorithm for finding the left-hand side.}
\label{finding lhs algorithm}
Suppose we have a word $x_1\cdots x_n \in \Astar$ and we know it has
a suffix which is the left-hand side of some rule in $\SetOfRules\Rules$.
Suppose no proper prefix of $x_1\cdots x_n$ has this property.
We give an algorithm that finds the shortest such suffix.

We read the word from right to left, starting with $x_n$.
We assume that $x_{k+1} x_{k+2} \cdots x_n$ has been read so far and that
as a result the current state of $Q(\Rules)$ is $S_k$, where
$S_k$ is a state of $Q(\Rules)$ (so $S_k$ is a subset of the set of states of
$\NRevLHS\Rules$).

We start the algorithm
with $k=n$ and the current state of $Q(\Rules)$ equal
to the singleton $\{(s_0,0,+)\}$ whose only element is the initial state
of $M$,
where $s_0$ is the initial state of $\Rev\Rules$.
$Q(\Rules)$ has three final states, namely the singleton sets
$\{f_i\}$ for $i=0,1,2$.

The steps of the algorithm are as follows:
\begin{enumerate}
\item\label{start} Record the current state as the $k$-th entry in an
array of size $n$, where $n$ is the length of the input word.
\item If the current state is not a final state, go to Step~\ref{jump}.
If the current state is a final state, then stop.
Note that the initial state of $Q(\Rules)$
is not a final state, so this step does not
apply at the beginning of the algorithm.
If the current state is a final state, then the
shortest suffix of $x_1\cdots x_n$ which is the left-hand side of a rule
in $\SetOfRules\Rules$ can then be proved to be $x_{k+1} x_{k+2} \cdots x_n$.
\item\label{jump}
If the arrow labelled $x_k$ with source the current state is
already defined, then redefine the current state to be the target of this arrow
and decrease $k$ by one.
\item If the preceding step does not apply, we have to compute the target
$T$ of the arrow labelled $x_k$ with source the current state $S_k$.
We do this by looking for all arrows labelled $x_k$
in $\NRevLHS\Rules$ with source in $S_k$. We define $T$ to be the set of
all targets of such arrows.
Note that this set of targets cannot be empty since we know
that some suffix of $x_1\cdots x_n$ is accepted by $\NRevLHS\Rules$.
\item\label{modified} There are two modifications which we can make to the
previous step.
\begin{enumerate}
\item\label{modifieda} Firstly, if the set of targets contains some final state
$f_j$, then we look for the largest value of $i=0,1,2$ such that $f_i\in T$
and redefine $T$ to be $\{f_i\}$.
We then insert into $Q(\Rules)$
an arrow labelled $x_k$ from $S_k$ to this final state.
If we have found that $T$ is a final state, we set $S_{k-1}$ equal to $T$,
decrease $k$ by one, and go to Step~\ref{start}.
\item\label{modifiedb}
Secondly, if, while calculating the set $T$, we find that a state $s$ of
$\Rev\Rules$ occurs in more than one triple $(s,i,j)$, then we only include
the triple with the largest value of $i$. For this to be well-defined,
we need to
know that $(s,i,+)$ and $(s,i,-)$ cannot both come up as potential
elements of
$T$---this is addressed in the proof of Theorem~\ref{lhs algorithm} along with
justifications of the other modifications.
\end{enumerate}
\item Having found $T$, see if it is equal to some state $T'$ of $Q(\Rules)$
which has already been constructed. If so, define an arrow labelled $x_k$
from $S$ to $T'$.
\item If $T$ has not already been constructed,
define a new state of $Q(\Rules)$
equal to $T$ and define an arrow labelled $x_k$ from $S$ to $T$.
\item
Set the current state equal to $T$ and decrease $k$ by one.
Then go to Step~\ref{start}.
\end{enumerate}

\begin{theorem}\label{lhs algorithm} Suppose $x_1\cdots x_n$ has a suffix
which is the left-hand side of a
rule in $\SetOfRules\Rules$ and suppose no prefix of $x_1\cdots x_n$
has this property.
Then the above algorithm correctly computes the shortest such
suffix.
\end{theorem}
\begin{proof}
We first show that the modification in Step~\ref{modifiedb} is well-defined in
the sense that triples $(s,i,+)$ and $(s,i,-)$ cannot both occur while
calculating $T$. The reason for this is that the third component can
only be $+$ if either none of $x_1\cdots x_n$ has been read, in which
case the only relevant state is $(s_0,0,+)$, or else only $x_n$
has been read, in which case the possible relevant states are $(f,1,-)$,
$(s,1,+)$ with $s\neq f$, and $(s,0,-)$. So a state of the form
$(s,i,j)$ with a given $s$ occurs at most once in a fixed subset with the
maximum possible value of $i$.

The effect of Step~\ref{modifieda} in the above algorithm is to ensure that
termination occurs as soon as a final state of $\Rev\Rules$ appears in a
calculated triple.
Since we know that $x_1\cdots x_n$ contains a left-hand
side of a rule in $\SetOfRules\Rules$ as a suffix we need only show that the
introduction of Step~\ref{modifiedb} does not affect the accepted language of
the constructed automaton. This will be a consequence of
Theorem~\ref{modified determinization}, as we now proceed to show.

Consider a triple $t=(s,i,j)$ arising during the calculation of a subset
$T$, and suppose that $s$ is a non-final state of $\Rev\Rules$.
If $j = +$ then $T$ cannot contain both $(s,0,+)$ and $(s,1,+)$ and so $t$ will
not be removed from $T$ as a result of Step~\ref{modifiedb}. Therefore we only
need to consider the case $j = -$. For $k=0,1,2$, let
$L_k \subseteq \Astar\times\Astar$ be the language
obtained by making $(s,k,-)$ the only initial state of $M$, and observe that
there can be no padded arrows in any path of arrows from $(s,k,-)$ to a
final state of $M$. Now by considering the definition of the non-padded
transitions in $M$ given in \ref{non-padded M transitions}, it is
straightforward to see that $L_0 \subseteq L_1 = L_2$.
Therefore, since $\NRevLHS\Rules$ has no $\epsilon$-arrows, we have just
shown that the hypotheses of Theorem~\ref{modified determinization} apply to
Step~\ref{modifiedb}. Hence the omission in
Step~\ref{modifiedb} does not affect the accepted language of $Q(\Rules)$.
\end{proof}

As with $P(\Rules)$, reading a word into $Q(\Rules)$ from right to
left can be slow in the initial stages of a Knuth--Bendix pass, but soon
speeds up to being linear with a small constant.

\subsection{Finding the right-hand side of a rule.}
\label{finding the right-hand side}
We retain the hypotheses of Section~\ref{properties of Rules}.
Namely, we have a two-variable rule automaton $\Rules$ which is welded
and satisfies various other minor conditions. We are given a word
$w=x_1\cdots x_n$, and we wish to reduce it relative to the rules implicitly
contained in $\Rules$.
So far we have located a left-hand side $\lambda$ which is
a subword of $w$. In this section we show how
to construct the corresponding right-hand side.

We first go into more detail as to how we propose to reduce $w$.
In outline we proceed as follows.

\subsection{Outline of the reduction process.}
\label{history stack}
\begin{enumerate}
\item Feed $w$ one symbol at a
time into the one-variable automaton $P(\Rules)$ described in
Section~\ref{Fast reduction},
storing the history of states reached on a stack.
\item If a final state is reached after some prefix $u$ of $w$ has been
read by $P(\Rules)$, then $u$ has some suffix which is a left-hand
side. Moreover, this procedure finds the shortest such prefix.
\item Feed $u$ from right to left into $Q(\Rules)$. A final state
is reached as soon as $Q(\Rules)$ has read the shortest suffix $\lambda$
of $u$ such that there is a rule $(\lambda,\rho) \in \SetOfRules\Rules$.
We now have $u = p\lambda$ and $w=p\lambda q$, where
$p,q \in \Astar$, every proper prefix of $p\lambda$
and every proper suffix of $\lambda$ is $\SetOfRules\Rules$-irreducible.
\item\label{use rule in Store}
Find $\rho$, the smallest word such that there
is a rule $(\lambda,\rho)$ in \Store (see \ref{Knuth--Bendix pass}).
If there is no such rule in \Store, find $\rho$ by a method to
be described in \ref{right-hand side routine},
such that $\rho$ is the smallest word
such that $(\lambda,\rho)\in \SetOfRules\Rules$.
\item\label{reduction gives new rule}
If $(\lambda,\rho)$ is not already in \Store, insert it into the
part of \Store called \New.
\item Replace $\lambda$ with $\rho$ in $w$ and pop $|\lambda|$ levels off
the stack so that the stack represents the history as it was immediately
after feeding
$p$ into $P(\Rules)$.
\item
Redefine $w$ to be $p\rho q$.
Restart at Step 1 as though $p$ has just been read and the next letter
to be read is the first letter of $\rho$. The history
stack enables one to do this.
\end{enumerate}

Note that other strategies might lead to
finding first some left-hand side in $w$ other than $\lambda$. Moreover,
there may be several different
right-hand sides $\rho$ with $(\lambda,\rho)\in
\SetOfRules\Rules$. A rule $(\lambda,\rho)$ in $\SetOfRules\Rules$ gives rise
to paths in $\Rules$, $\SLtwo$ and $\DRev\SLtwo$. We will find the path for
which right-hand side $\rho$ is \shl-least, given
that the left-hand side is equal to $\lambda$.

Let $\lambda = y_1\cdots y_{m}$.
Recall that a state of the one-variable automaton $Q(\Rules)$ used to find
$\lambda$ is a set of states of the form
$(s,i,j)$, where $s$ is a state of $\Rules, i\in\{0,1,2\}$ and $j\in\{+,-\}$.
When finding $\lambda$ we kept the history of states of $Q(\Rules)$ which were
visited---see Step~\ref{start}.
Let $Q_k$ be the set of triples $(s,i,j)$
comprising the state of $Q(\Rules)$ after reading
the word $ y_{k+1}\cdots y_m$ from right to left.
$Q_0=\{f_i\} = \{(s_0,i,-)\}$ where $s_0$
is the unique initial and final state of $\Rules$, and $i$ is the
difference in length between $\lambda$ and the $\rho$ that we are looking for.

\subsection{Right-hand side routine.}
\label{right-hand side routine}
Inductively, after reading $y_1\cdots y_k$ we will have determined
$z_1\cdots z_k$, the prefix of $\rho$. Inductively we also have a triple
$(s_k,i_k,j_k)$, where $s$ is a state of $\Rules$, $i_k$ is 0 or 1 or
2 and $j_k$ is $+$ or $-$. Note that we always have $m-k \ge i_k$.
\begin{enumerate}
\item\label{start2} If $m-k = i_k$, then we have found
$\rho = z_1\cdots z_k$ and we stop. So from now on we assume that
$m > i_k + k$. This means that the next symbol $(y_{k+1},z_{k+1})$ of
$(\lambda,\rho)$ does not have a padding symbol in its right-hand component.
\item We now try to find $z_{k+1}$ by running through each element $z\in A$ in
increasing order. Set $z$ equal to the least element of $A$.
\item \label{possible equality}
If $k=0$ and $i_0=0$, then $\lambda$ and $\rho$ will be of
equal length, so the first symbol of $(\lambda,\rho)$ must be
$(y_1,z_1)$, where $y_1> z_1$.
So at this stage we can prove that we have $y_1 > z$, since we know
that there must be some right-hand side corresponding to our given
left-hand side.

If $k=0$ and $i_0 > 0$, then the first symbol of
$(\lambda,\rho)^+$ is $(y_1,z_1)$ with $z_1\in A$ and $y_1\neq z_1$.
If $k=0$, $i_0>0$ and $y_1=z$, we increase $z$ to the next element of $A$.
\item\label{find z}
Here we are trying out a particular value of $z$
to see whether it allows us to get further.
We look in $\Rules$ to see if $s_k^{(y_{k+1},z)} = s_{k+1}$ is defined.
If it is not defined, we increase $z$ to the next element of $A$
and go to Step~\ref{possible equality}.

\item\label{finding the triple}
If $s_{k+1}$ is defined in Step~\ref{find z}, we look in $Q_{k+1}$ for a
triple $(s_{k+1},i_{k+1},j_{k+1})$ which is the source of an arrow
labelled $(y_{k+1},z)$ in the automaton $M$, defined in
Section~\ref{finding the left-hand side}.
Note that, by the proof of \thref{lhs algorithm}, $Q_{k+1}$ contains at most
one element whose first coordinate is $s_{k+1}$. As a result, the search can be
quick.
\item If $(s_{k+1},i_{k+1},j_{k+1})$ is not found in Step~\ref{finding the triple},
increase $z$ to the next element of $A$ and go to Step~\ref{possible equality}.
\item If $(s_{k+1},i_{k+1},j_{k+1})$ is found in Step~\ref{finding the triple},
set $z_{k+1} = z$, increase $k$ and go to Step~\ref{start2}.
\end{enumerate}

The above algorithm will not hang, because each triple $(s_k,i_k,j_k)$
that we use does come from a path of arrows in $M$ which starts
at the initial state of $M$ and ends at the first possible final state
of $M$. Therefore all possible right-hand sides $\rho$ such that
$(\lambda,\rho) \in \SetOfRules\Rules$, are implicitly computed when we
record the states of $Q(\Rules)$ (see
Step~\ref{start}).
Since $i_k$ does not vary during our search, we
will always find the shortest possible $\rho$, with $|\lambda| -
|\rho| $ being equal to this constant value of $i_k$.
Since we always look for $z$ in increasing order, we are
bound to find the lexicographically least $\rho$.

%% file: Sec8.tex
\section{A modified determinization algorithm}
\label{A modified determinization algorithm}[Section]

In this section we discuss a useful modification to the usual determinization
algorithm for turning an NFA into a DFA.
Let $N$ be an NFA. The usual proof that $N$ can be determinized,
is to form a new automaton $M$ each state of which is a subset $\sigma$
of the set $S(N)$
of states of $N$ such that $\sigma$ is $\epsilon$-closed.
That is to say,
if $s \in \sigma \subset S(N)$, then each $\epsilon$-arrow with source
$s$ also has target in $\sigma$. The initial state of $M$
is the $\epsilon$-closure
of the set of all initial states in $N$. The effect of an arrow labelled
$x \in A$ on $\sigma$ is to take each $s\in\sigma$, apply $x$ in all
possible ways, and then to take the $\epsilon$-closure of the subset
of $S(N)$ so obtained. A final state of $M$ is any subset of $S(N)$
containing a final state of $N$.

In practice, to find $M$,
we start with the $\epsilon$-closure of the set of initial states of
$N$ and proceed inductively. If we have found a state $s$ of $M$ as a
subset of the set of states of $N$, we fix some $x\in A$, and
apply $x$ in all possible ways to all $t\in s$, where $t$ is a state
of $N$. We then follow with $\epsilon$-arrows to form an $\epsilon$-closed
subset of states of $N$. This gives us the result of applying $x$ to
$s$. The modification we wish to make to the usual subset
construction is now explained and justified.

We will denote by $M'$ the modified version of $M$ thus obtained. $M'$ is
a DFA which accepts the same language as $M$ and $N$, but the structure of
$M'$ might be simpler than that of $M$.

Suppose $p$ is a state of the NFA $N$. Let $N_p$ be the same automaton
as $N$, except that the only initial state is $p$. Suppose $p$ and $q$
are distinct states of $N$ and that $L(N_p)\subset L(N_q)$. Suppose
also that the $\epsilon$-closure of $q$ does not include $p$. Under
these circumstances, we can modify the subset construction as follows.
As before, we start with the $\epsilon$-closure of the set of initial
states of $N$. We follow the same procedure for defining the arrows and
states of $M'$ as for $M$, except that, whenever we construct a subset
containing both $p$ and $q$, we change the subset by omitting $p$.

\subsection{Required conditions.}\label{required conditions}
The situation can be generalized. We suppose that we have a partial order
defined on the set of states of $N$, such that, if $p < q$, then
$L(N_p) \subset L(N_q)$.
We assume that if $p < q$, $p'<q'$ and $p'$ is contained in
the $\epsilon$-closure of $q$, then $p'=q$.

We follow the same procedure for
defining the arrows and states of $M'$ as for $M$, except that, whenever
we construct a subset containing both $p$ and $q$ with $p<q$, we change the
subset by omitting $p$.

\begin{theorem}\label{modified determinization}
Under the above hypotheses, $L(M') = L(N)$.
\end{theorem}
\begin{proof}
Consider a word $w=x_1\cdots x_n\in \Astar$ which is accepted by $N$
via the path of arrows in $N$
$$(v_0,\epsilon\ast,u_1,x_1,v_1,\cdots,v_{n-1},\epsilon\ast,u_n,x_n,v_n,
\epsilon\ast,u_{n+1}).$$
This means that,
for each $i$ with $0\le i \le n$,
there is an $x_i$-arrow in $N$ from
$u_i$ to $v_i$ and
$u_{i+1}$ is in the $\epsilon$-closure of $v_i$.
Moreover $v_0$ is an initial state and $u_{n+1}$ is a final state.

Our proof will be by induction on $i$.
The $i$-th statement in the induction is that
we have states $s_0,\ldots,s_i$ of
$M'$ such that $s_0$ is the initial state and, for each $j$ with
$0<j < i$, there
is an arrow $x_j:s_{j-1} \to s_j$ in $M'$, so that,
after reading $x_1\cdots x_{i-1}$,
$M'$ is in state $s_{i-1}$.
Our induction statement also says that we have
a path of arrows in $N$
$$(u_i^i,x_i,v_i^i,\epsilon\ast,u_{i+1}^i,\cdots,
u_n^i,x_n,v_n^i,\epsilon\ast,u_{n+1}^i),$$
such that $u_i^i\in s_{i-1}$ and $u_{n+1}^i$ is a final state of $N$.

The induction starts with $i=1$ and $s_0$ the initial state of $M'$.
We form $s_0$ by taking all initial states of $N$, and taking their
$\epsilon$-closure. If this subset of states of $N$ contains both
$p$ and $q$ with $p<q$, then $p$ is omitted from $s_0$, the initial state
of $M'$.
If $u_1\notin s_0$, then we must have $u_1 = p$, with $q\in s_0$ and
$p<q$.
So $q$ must be a maximal element of $s_0$ with respect to the
partial order.
Now $w\in L(N_{p}) \subset L(N_{q})$.
It follows that
we can take $u_1^1$ in the $\epsilon$-closure of $q$ and then define the
rest of the path of arrows for the case $i=1$.
Since $q\in s_0$ and
$u_1^1$ is in the $\epsilon$-closure of $q$, it is not the case that
there is a $q'$ such that
$u_1^1 < q'\in s_{0}$, according to \ref{required conditions}.
So $u_1^1\in s_0$ (that is, it is not omitted in our construction)
and the induction can start.

Now suppose the induction statement is true for $i$.
We prove it for $i+1$.
we have a path of arrows
$$(u_i^i,x_i,v_i^i,\epsilon\ast,u_{i+1}^i,\cdots,
u_n^i,x_n,v_n^i,\epsilon\ast,u_{n+1}^i),$$
in $N$ such that $u_i^i\in s_{i-1}$ and $u_{n+1}^i$ is a final state of $N$.
We define $s_i$ from $s_{i-1}$ in the manner described above.
First we
apply $x_i$ in all possible ways to all states in $s_{i-1}$,
obtaining $v_i^i$ as one of the target states,
and then take the $\epsilon$-closure,
obtaining $u_{i+1}^i$ as one of the targets of an $\epsilon$-arrow.
Finally, if
$s_i$ contains both $p$ and $q$, with $p < q$ then $p$ is deleted
from $s_i$ before $s_{i}$ becomes a state of $M'$.

It now follows that either $u_{i+1}^i\in s_i$, or else, for some $p<q$,
$u_{i+1}^i = p$, $q\in s_i$ and $p\notin s_i$.
In the first case we define $u_j^{i+1} = u_j^i$ and $v_j^{i+1} = v_j^i$ for
$j> i$ and the induction step is complete.
In the second case, using the fact that $x_{i+1}\cdots x_n \in L(N_{p})
\subset L(N_q)$, we see that
we can take $u_{i+1}^{i+1}$ in the $\epsilon$-closure of $q$
and then define the rest of the path of arrows.
Since $q\in s_i$ and $u_{i+1}^{i+1}$ is in the $\epsilon$-closure of $q$,
\ref{required conditions} shows that it is not possible to have
$q'\in s_i$ and $u_{i+1}^{i+1}<q'$.
Therefore $ u_{i+1}^{i+1}\in s_i$. This completes the induction step.

At the end of the induction, $M'$ has read all of $w$ and is in state $s_n$.
We also have the final state $u_{n+1}^{n+1} \in s_n$, so that $w$ is accepted
by $M'$.

Conversely, suppose $w$ is accepted by $M'$. It follows easily by induction
that if $M'$ is in state $s_i$ after reading the prefix $x_1\cdots x_i$ of
$w$, then each state $u\in s_i$ can be reached from some initial state of
$N$ by a sequence of arrows labelled successively
$x_1,\ldots,x_i$, possibly interspersed with $\epsilon$-arrows.
Now $s_n$ must contain a final state, and so $w$ is accepted by $N$.
\end{proof}

\begin{remark}
The practical usage of this theorem clearly depends on having an
efficient way of determining when the condition $L(N_p) \subset L(N_q)$
is satisfied. In this paper we have seen several
examples of such tests which cost
virtually nothing to implement but have the potential to save an appreciable
amount of both space and time.
\end{remark}

%% file: Sec9.tex
\section{Miscellaneous details}
In this section we present a number of points which did not seem to
fit elsewhere in this paper.

\subsection{Aborting.} It is possible that we come to a situation where
the procedure is not noticing that certain words are reducible, even though
the necessary information to show that they are reducible is already
in some sense known. It is also possible that reduction is being carried
out inefficiently, with several steps being necessary, whereas in some
sense the necessary information to do the reduction in one step is
already known. An indication that our procedure is not proceeding as
well as one hoped might be that
\WDiff is constantly changing, with states being identified and consequent 
welding, or with new states or arrows being added. In this case it might be 
advisable to abort the current Knuth--Bendix pass.

To see if abortion is advisable, we can record statistics about how much
\WDiff has changed since the beginning of a pass. If the changes seem
excessive, then the pass is aborted. A convenient place for the program to
decide to do this is just before another rule from \New is examined
at Step~\ref{process New}.

If an abort is decided upon then all
states and arrows of \WDiff are marked as \needed.
At this point the program jumps to Step~\ref {start main loop}.

\subsection{Priority rules.} A well-known phenomenon found when
using Knuth--Bendix to look for automatic structures, is that rules
associated with finding new word differences or new arrows in
\WDiff should be used more intensively than other rules. Further aspects
of the structure are then found more quickly.
This is not a theorem---it is observed behaviour seen on examples
which happen to have been investigated.

A new rule associated with new word differences or new arrows in
\WDiff is
marked as a priority rule. When a priority rule is minimized, the output is
also marked as a priority rule. If a priority rule is added to one of the lists
\Considered, \This or \New, it is added to the front of the list, whereas
rules are normally added to the end of the list.
Just before deciding to add a priority rule to \New, we check to see if
the rule is minimal. If so, we add it to the front of \This instead of to
the front of \New.

When a rule is taken from \This at Step~\ref{process This} during the 
main loop, it is normally compared with all rules in \Considered,
looking for overlaps between left-hand sides. In the case of a priority rule,
we compare left-hand sides not only with rules in \Considered,
but also with all rules in \This. If a normal rule $(\lambda,\rho)$
is taken from \This and comparison with a rule in \Considered gives rise
to a priority rule, then the rule $(\lambda,\rho)$ is also marked as a priority
rule. It is then compared with all rules in \This, once it has been
compared with all rules in \Considered.

Treating some rules as priority rules makes little difference unless there
is a mechanism in place for aborting a Knuth--Bendix pass when \WDiff
has sufficiently changed. If there is such a mechanism, it can make a big
difference.

\subsection{An efficiency consideration.}
During reduction we often have a state $s$ in a two-variable
automaton and an $x\in A$, and we are looking for an arrow
labelled $(x,y)$ with certain properties, where $y\in A^+$. It therefore
makes a big difference if the arrows with source $s$ are arranged so that
we have rapid access to arrows labelled $(x,y)$ once $x$ is given.

\subsection{The present.}
Many of the ideas in this paper have been implemented
in C++ by the second author.
But some of the ideas in this paper only occurred to
us while the paper was being written, and the procedures and algorithms
presented in this paper seem to us to be substantial improvements on what
has been implemented so far. An unfortunate result of this is that we are
unable to present experimental data to back up our ideas, although many
of our ideas have been explored in depth with actual code. Our experimental 
work has been essential in enabling us to come to the
better algorithms which are presented here.

\subsection{Comparison with \textit{kbmag}.} Here we describe the
differences between our ideas and the ideas in Derek Holt's \textit{kbmag}
programs \cite{Holt:KBMAG}. These programs try to compute the 
\shl-automatic
structure on a group. Our program is a substitute only for the first program
in the \textit{kbmag} suite of programs.

In \textit{kbmag}, fast reduction is carried out using an automaton with
a state for every prefix of every left-hand side. In our program we also
keep every rule. However, the space required by a single character in our
program is less by a constant multiple than the space required for a
state in a finite state automaton. Moreover, compression techniques
could be used in our situation so that less space is used, whereas
compression is not available in the situation of \textit{kbmag}.

The other large objects in our set-up are the automata $P(\Rules[n])$
defined in \ref{P} and $Q(\Rules[n])$ defined in \ref{Q}.
In \textit{kbmag}, there has also to be an automaton like $P(\Rules[n])$, and
it is possible to arrange that this automaton is only constructed
after the Knuth--Bendix process is halted.
In \textit{kbmag} there is no analogue of our $Q(\Rules[n])$.
So these are advantages of \textit{kbmag}.

In \textit{kbmag}, reduction is carried out extremely rapidly. However,
as new rules are found, the automaton in \textit{kbmag} needs to be
updated, and this is quite time-consuming. In our situation, updating the
automata is quick, but reduction is slower by a factor of around three,
because the word has to be
read into two or three different automata. Moreover we sometimes need to use the
method of Section~\ref{finding the right-hand side} which is slower
(by a constant factor) than simply reading a word into a deterministic
finite state automaton.

In \textit{kbmag}, there is a heuristic, which seems to be inevitably arbitrary,
for deciding when to
stop the Knuth--Bendix process. In our situation there is a sensible
heuristic, namely we stop if we find $\Rules[n+1] = \Rules[n]$.

In the case of \textit{kbmag}, there are
occasional cases where the process of finding the set of word differences
oscillates indefinitely. This is because redundant rules are sometimes
unavoidably introduced into the set of rules, introducing unnecessary
word differences. Later redundant rules are eliminated and also the
corresponding word differences. This oscillation can continue indefinitely.
Holt has tackled this problem in his programs
by giving the user interactive modes of running them.

In our case,
 the results in Section~\ref{correctness} show that, given a
\shl-automatic group, the automaton $\Rules[n]$ will eventually stabilize,
as proved in \thref{main},
given enough time and space.

We believe that the main advantage of our approach for computing
automatic structures will only become
evident (if it exists at all)
when looking at very large examples. We plan to carry out a
systematic examination of \shl-automatic groups generated
by Jeff Weeks' \textit{SnapPea} program---see
\cite{Weeks:Snappea}---in order to carry out a systematic comparison.